\documentclass[oneside,english]{amsart}
\usepackage[T1]{fontenc}
\usepackage[latin9]{inputenc}
\usepackage{amsthm}
\usepackage{amstext}
\usepackage{amssymb}
\usepackage{setspace}
\usepackage{esint}
\makeatletter

\newcommand{\noun}[1]{\textsc{#1}}

\numberwithin{equation}{section}
\numberwithin{figure}{section}
 \theoremstyle{definition}
 \newtheorem*{defn*}{\protect\definitionname}
\theoremstyle{plain}
\newtheorem{thm}{\protect\theoremname}
  \theoremstyle{plain}
  \newtheorem{lem}[thm]{\protect\lemmaname}
  \theoremstyle{remark}
  \newtheorem*{rem*}{\protect\remarkname}
  \theoremstyle{plain}
  \newtheorem{cor}[thm]{\protect\corollaryname}
  \theoremstyle{plain}
  \newtheorem*{lem*}{\protect\lemmaname}

\date{}

\makeatother

\usepackage{babel}
  \providecommand{\corollaryname}{Corollary}
  \providecommand{\definitionname}{Definition}
  \providecommand{\lemmaname}{Lemma}
  \providecommand{\remarkname}{Remark}
\providecommand{\theoremname}{Theorem}

\begin{document}

\title{Compound Poisson statistics in conventional and nonconventional setups}

\author{Ariel Rapaport}

\date{April 14, 2014}

\subjclass[2000]{Primary: 60F05 Secondary: 37D35, 60J05}

\keywords{$\psi$-mixing, Pólya\textendash{}Aeppli distribution, Gibbs measure,
$f$-expansion}

\thanks{Supported by ERC grant 306494}
\begin{abstract}
Given a periodic point $\omega$ in a $\psi$-mixing shift with countable
alphabet, the sequence $\{S_{n}\}$ of random variables counting the
number of multiple returns to shrinking cylindrical neighborhoods
of $\omega$ is considered. Necessary and sufficient conditions for
the convergence in distribution of $\{S_{n}\}$ are obtained, and
it is shown that the limit is a Pólya\textendash{}Aeppli distribution.
A global condition on the shift system, which guarantees the convergence
in distribution of $\{S_{n}\}$ for every periodic point, is introduced.
This condition is used to derive results for $f$-expansions and Gibbs
measures. Results are also obtained concerning the possible limit
distribution of sub-sequences $\{S_{n_{k}}\}$. A family of examples
in which there is no convergence is presented. We exhibit also an
example for which the limit distribution is pure Poissonian.
\end{abstract}
\maketitle

\section{Introduction}

In this article $(\Omega,T)$ is a shift space with a countable alphabet,
equipped with a $\psi$-mixing and $T$-invariant probability measure
$\mathbb{P}$. Given $\omega\in\Omega$, the sequence of random variables
$\{S_{n}^{\omega}\}_{n=1}^{\infty}$ is considered, where for each
$n\geq1$, 
\[
S_{n}^{\omega}(\gamma)=\sum_{k=1}^{N_{n}^{\omega}}\prod_{j=1}^{\ell}1_{A_{n}^{\omega}}\circ T^{d_{j}k}(\gamma)\quad\mbox{ for }\gamma\in\Omega,
\]
$A_{n}^{\omega}$ is the cylinder set determined by the $n$-prefix
of $\omega$, $N_{n}^{\omega}=[(\mathbb{P}(A_{n}^{\omega}))^{-\ell}]$,
and $1\leq d_{1}<...<d_{\ell}$ are fixed integers. When the setup
is said to be conventional it means that $\ell=1$ and $d_{1}=1$.
We are interested in the limit distribution $\mu_{\omega}$, of the
sequence $\{S_{n}^{\omega}\}$.

In the conventional case, the existence and characterization of $\mu_{\omega}$
is a widely studied problem. Poisson approximation estimates for almost
all non-periodic $\omega$, were obtained in \cite{key-14} and \cite{key-15}.
For a periodic $\omega$, compound Poisson approximations were derived
in \cite{key-13}. In \cite{key-16}, it was shown that the distributional
limit of the normalized number of returns to small neighborhoods of
periodic points is compound Poisson. This was shown for certain non-uniformly
hyperbolic dynamical systems, which include certain piecewise expanding
maps of the interval.

The problem, in the nonconventional setup described above, was considered
for the first time in \cite{key-17}. It was shown there that if $(\Omega,T)$
is a subshift of finite type and $\mathbb{P}$ is a Gibbs invariant
measure, then $\mu_{\omega}$ exists (i.e. $\{S_{n}^{\omega}\}$ has
a limit in distribution) and is equal to the Poisson distribution
with the parameter $1$, for $\mathbb{P}$-almost all $\omega\in\Omega$.
In \cite{key-1} this result was extended for every non-periodic $\omega$
under more general $\psi$-mixing assumptions. It was also shown there
that if $\omega$ is a periodic point and the dynamical system is
a Bernoulli shift with a countable state space, then $\mu_{\omega}$
is a Pólya\textendash{}Aeppli distribution.

Here the case of a periodic $\omega$ is considered for a general
$\psi$-mixing system in the nonconventional setup. It is shown that
if $\mu_{\omega}$ exists, then it must be a Pólya\textendash{}Aeppli
distribution (which may be pure Poissonian). A sufficient condition
for this is given, which was first introduced for the conventional
setup in \cite{key-13}. A necessary condition is also given which
is the same as the sufficient condition for the conventional setup.
Both the necessary condition and the fact that $\mu_{\omega}$ must
be a Pólya\textendash{}Aeppli distribution (if it exists) are new
also for the conventional case.

A verifiable global condition on the dynamical system that guarantees
the existence of $\mu_{\omega}$ for each periodic $\omega$ is introduced.
This condition is based on the notion of the inverse Jacobian of a
measure preserving transformation, which is defined in the next section.
Using this, it is shown that if the shift system is derived from an
$f$-expansion on $[0,1]$ with its absolutely continuous invariant
measure (see \cite{key-2}) then $\mu_{\omega}$ always exists, in
particular this applies to the Gauss map equipped with the Gauss measure.
For the conventioanl case this follows from results found in \cite{key-16}.
The global condition is also used for showing that if $(\Omega,T)$
is of finite type and $\mathbb{P}$ is a Gibbs measure then $\mu_{\omega}$
always exists. This was first shown for the conventional setup in
\cite{key-13}.

If $\omega$ is periodic then the sequence $\{S_{n}^{\omega}\}$ does
not necessarily have a limit in distribution, an example of this phenomena
was given in \cite{key-1}. However, the partial limits in distribution
of $\{S_{n}^{\omega}\}$ can be characterized. It is shown here that
if a sub-sequence $\{S_{n_{k}}^{\omega}\}$ converges in distribution,
then the limit distribution must be a compound Poisson distribution.
Also, the nonconvergence example from \cite{key-1} is extended into
a wide class of examples in the same spirit. In addition, an example
of a finite type system and a periodic point $\omega$ is given, for
which $\mu_{\omega}$ equals a pure Poisson distribution.

The rest of this article is organized as follows: In Section \ref{S2}
the notation and framework being used are presented. In Section \ref{S3}
the results are stated. In Section \ref{S4} we prove the results
regarding the pointwise conditions. In Section \ref{S5} we prove
the result concerning the global condition, and the two applications
of it are obtained. In Section \ref{S6} the results regarding the
classification of partial limits are proved. In Section \ref{S7}
the family of nonconvergence examples is developed. In Section \ref{S8}
we construct the example of a periodic point for which the limit distribution
is pure Poissonian.

\textbf{Acknowledgment. }I would like to thank my adviser Professor
Yuri Kifer, for suggesting to me problems studied in this paper, and
for many helpful discussions.

\section{\label{S2}Framework and notations}

\subsection{\label{SS2}The underlying dynamical system}

Let $\mathcal{A}$ be a finite or countable set (the alphabet), with
$|\mathcal{A}|>1$. Let $\left(\Omega,\mathcal{F},\mathbb{P},\mathbb{\mathit{T}}\right)$
be a measure-preserving system where $\Omega=\mathcal{A}^{\mathbb{N}}$,
$\mathcal{F}$ is the $\sigma$-algebra generated by the coordinates
projections from $\Omega$ onto $\mathcal{A}$, $\mathbb{\mathit{T}}$:$\Omega\rightarrow\Omega$
is the left shift, and $\mathbb{P}$ is a $\mathbb{\mathit{T}}$-invariant
probability measure. $\mathcal{A}$ is considered as a topological
space with the discrete topology, and $\mathcal{A}^{\mathbb{N}}$
as a topological space with the product topology. For distinct $\omega,\gamma\in\Omega$
define $d(\omega,\omega)=0$ and $d(\omega,\gamma)=2^{-m}$, where
$m=\min\{j\geq0\::\:\omega_{j}\ne\gamma_{j}\}$. Then $d$ is a metric
on $\Omega$ which induces the product topology.

Given $J\subset\mathbb{N}$ define
\[
\mathcal{F}_{J}=\sigma\{\{\omega_{j}=a\}\::\: j\in J\mbox{ and }a\in\mathcal{A}\}
\]
We assume that $\mathbb{P}$ is $\psi$-mixing, i.e. that there exists
a sequence $\left\{ \psi_{m}\right\} _{m\geq0}\subset\mathbb{R}^{+}$
with $\psi_{m}\overset{m\rightarrow\text{\ensuremath{\infty}}}{\longrightarrow}0$,
such that for each $m,n\in\mathbb{N}$, $E\in\mathcal{F}_{\{0,...,n-1\}}$
and $F\in\mathcal{F}$,
\begin{equation}
|\mathbb{P}(E\cap T^{-(n+m)}F)-\mathbb{P}(E)\mathbb{P}(F)|\leq\psi_{m}\mathbb{P}(E)\mathbb{P}(F)\label{E18}
\end{equation}
Examples of $T$-invariant and $\psi$-mixing measures include Gibbs
measures on a subshift of finite type (see \cite{key-12}), and measures
derived from $f$-expansions defined on the unit interval (see \cite{key-2}).

Sometimes it will be more convenient to work in a certain closed subspace
of $\Omega$, namely in a topological Markov shift (see \cite{key-10}),
which will now be defined. Let $S=(S_{a,b})_{a,b\in\mathcal{A}}$
be a matrix of $0$'s and $1$'s with no columns or rows which are
all $0$'s. Let
\[
\Omega_{S}=\{\omega\in\Omega\::\: S_{\omega_{j},\omega_{j+1}}=1\mbox{ for each }j\geq0\}
\]
and 
\[
\mathcal{F}_{S}=\{E\in\mathcal{F}\::\: E\subset\Omega_{S}\}
\]
Then $\Omega_{S}$ is a closed and $T$-invariant subset of $\Omega$,
which is considered as a topological space with the subspace topology
inherited from $\Omega$. It is assumed that $\mathbb{P}(\Omega_{S})=1$,
and that $T:\Omega_{S}\rightarrow\Omega_{S}$ is topologically mixing,
i.e. for every pair of open sets $U,V\subset\Omega_{S}$ there exist
a number $N(U,V)\in\mathbb{N}_{+}$ such that $U\cap T^{-n}V\ne\emptyset$
for all $n\geq N(U,V)$.

\subsection{Words and cylinders}

Let $\mathcal{A}^{*}$ be the set of finite words over $\mathcal{A}$,
and let $\mathcal{A}_{S}^{*}\subset\mathcal{A}^{*}$ be the subset
of $S$-admissible words, i.e.
\[
\mathcal{A}_{S}^{*}=\{a_{0}\cdot...\cdot a_{r-1}\in\mathcal{A}^{*}\::\: S_{a_{j-1},a_{j}}=1\mbox{ for each }1\leq j<r\}
\]
For each $u,w\in\mathcal{A}^{*}$, $\omega\in\Omega$ and $k\geq0$,
let $u\cdot w\in\mathcal{A}^{*}$ be the concatenation of $u$ and
$w$, let $w^{k}\in\mathcal{A}^{*}$ be the concatenation of $w$
with itself $k$ times, and let $u\cdot\omega\in\Omega$ be the sequence
obtained by adding $u$ to the beginning of $\omega$. Given $a_{0},....,a_{r-1}\in\mathcal{A}$,
$a_{0}\cdot....\cdot a_{r-1}=w\in\mathcal{A}^{*}$ and $n\geq1$ let
\[
w^{n/r}=w^{[n/r]}\cdot a_{0}\cdot...\cdot a_{n-r[\frac{n}{r}]-1}
\]
where $[\frac{n}{r}]$ stands for the integral part of $\frac{n}{r}$,
and let
\[
[w]=\{\omega\in\Omega\::\:\omega_{j}=a_{j}\mbox{ for each }0\leq j<r\}.
\]
The set $[w]$ is called an $r$-cylinder. Sometimes $[w]^{n/r}$
is written in place of $[w^{n/r}]$, and $[a_{0},...,a_{r-1}]$ in
place of $[w]$.

From the $\psi$-mixing assumption and Lemma 3.1 in \cite{key-1}
it follows that there exist a constant $\Gamma>0$ such that 
\begin{equation}
\mathbb{P}[a_{0},...,a_{n-1}]\leq e^{-\Gamma n}\label{E2}
\end{equation}
for each $a_{0},...,a_{n-1}\in\mathcal{A}$.

Given an $r$-cylinder $A$ let $\pi(A)$ be the period of $A$, i.e.
\[
\pi(A)=\min\{j\in\{1,...,r\}\::\: A\cap T^{-j}A\ne\emptyset\}
\]
Given $\omega\in\Omega$ and $n\geq1$, define $A_{n}^{\omega}=[\omega_{0},...,\omega_{n-1}]$.
Let $\Omega_{\mathbb{P}}\subset\Omega_{S}$ be the support of $\mathbb{P}$,
then
\[
\Omega_{\mathbb{P}}=\{\omega\in\Omega\::\:\mathbb{P}(A_{n}^{\omega})>0\mbox{ for each }n\geq1\}
\]

\subsection{The observables}

The random variables counting the number of multiple recurrences to
cylindrical neighborhoods will now be defined. Let $1\leq d_{1}<...<d_{\ell}$
be integers. For each cylinder $A\in\mathcal{F}$, $N\in\mathbb{N}_{+}:=\{1,2,...\}$
and $\omega\in\Omega$ set
\[
S_{N}^{A}(\omega)=\sum_{k=1}^{N}X_{k}^{A}(\omega)\mbox{ where }X_{k}^{A}(\omega)=\prod_{j=1}^{\ell}1_{A}\circ T^{d_{j}k}(\omega)\mbox{ for each }k\in\mathbb{N}_{+}
\]
where $1_{A}$ stands for the indicator function of the set $A$.
As mentioned above, when we say that the setup is conventional it
means that $\ell=1$ and $d_{1}=1$.

\subsection{The inverse Jacobian}

For each measurable $E\in\mathcal{F}_{S}$ define
\[
\mathbb{P}\circ T(E)=\sum_{a\in\mathcal{A}}\mathbb{P}(T(E\cap[a]))
\]
then $\mathbb{P}\circ T$ is a $\sigma$-finite measure on $(\Omega_{S},\mathcal{F}_{S})$,
which is finite on cylinders (see \cite{key-10} for more details
on $\mathbb{P}\circ T$). Given $E\in\mathcal{F}_{S}$,
\[
\mathbb{P}\circ T(E)=\sum_{a\in\mathcal{A}}\mathbb{P}(T^{-1}(T(E\cap[a])))\geq\sum_{a\in\mathcal{A}}\mathbb{P}(E\cap[a])=\mathbb{P}(E),
\]
so $\mathbb{P}\ll\mathbb{P}\circ T$, and the following definition
makes sense.
\begin{defn*}
The function $J=\frac{d\mathbb{P}}{d\mathbb{P}\circ T}\in L^{1}(\Omega_{S},\mathcal{F}_{S},P\circ T)$
will be called here the inverse Jacobian of $\mathbb{P}$. 
\end{defn*}

The global condition mentioned in the introduction involves the concept
of the inverse Jacobian. More on this notion can be found in \cite{key-10}
and \cite{key-11}. In \cite{key-10} the function $\frac{d\mathbb{P}}{d\mathbb{P}\circ T}$
is called the Jacobian, and in \cite{key-11} this name is given to
$\frac{d\mathbb{P}\circ T}{d\mathbb{P}}$. Here we use the function
$\frac{d\mathbb{P}}{d\mathbb{P}\circ T}$, since always $\mathbb{P}\ll\mathbb{P}\circ T$
whenever $\mathbb{P}$ is $T$-invariant.

\subsection{Probability measures on $\mathbb{N}$}

Let $\mathcal{M}(\mathbb{N})$ denote the collection of all probability
measures on $\mathbb{N}$. Given a random variable $Y$, the distribution
of $Y$ is denoted by $\mathcal{L}(Y)$. Given $\mu\in\mathcal{M}(\mathbb{N})$,
it is written $Y\sim\mu$ if $Y$ is a random variable with $\mathcal{L}(Y)=\mu$.
Given random variables $Y,Y_{1},Y_{2},...$ we write $Y_{n}\overset{d}{\Longrightarrow}Y$
if the sequence $\{Y_{j}\}_{j=1}^{\infty}$ converges to $Y$ in distribution.

The total variation distance between members of $\mathcal{M}(\mathbb{N})$
is denoted by $d_{TV}$, i.e. given $\mu,\nu\in\mathcal{M}(\mathbb{N})$,
\[
d_{TV}(\mu,\nu)=\sup\{|\mu(E)-\nu(E)|\::\: E\subset\mathbb{N}\}
\]
Given $\mu,\mu_{1},\mu_{2},...\in\mathcal{M}(\mathbb{N})$ we write
$\mu_{j}\overset{d}{\Longrightarrow}\mu$ if the sequence $\{\mu_{j}\}_{j=1}^{\infty}$
converges to $\mu$ in distribution. Then $\mu_{j}\overset{d}{\Longrightarrow}\mu$
if and only if $d_{TV}(\mu_{j},\mu)\overset{j}{\rightarrow}0$, which
holds if and only if $|\mu\{k\}-\mu_{j}\{k\}|\overset{j}{\rightarrow}0$
for each $k\in\mathbb{N}$. See \cite{key-4} for more details on
the total variation distance.

A sequence $\{\mu_{j}\}_{j=1}^{\infty}\subset\mathcal{M}(\mathbb{N})$
is said to be tight if for every $\epsilon>0$ there exists $N\geq1$
such that $\mu_{j}[N,\infty)\leq\epsilon$ for each $j\geq1$. It
holds that $\{\mu_{j}\}_{j=1}^{\infty}$ is tight if and only if for
every sub-sequence $\{\mu_{j_{k}}\}_{k=1}^{\infty}$ there exist a
further sub-sequence $\{\mu_{j_{k_{i}}}\}_{i=1}^{\infty}$ and $\mu\in\mathcal{M}(\mathbb{N})$,
such that $\mu_{j_{k_{i}}}\overset{d}{\Longrightarrow}\mu$ as $i\rightarrow\infty$
(see Theorem 25.10 in \cite{key-9}).

Given a random variable $Y$, the characteristic function of $Y$
is denoted by $\varphi_{Y}$, i.e. $\varphi_{Y}(x)=E[e^{ixY}]$ for
each $x\in\mathbb{R}$. Given $\mu\in\mathcal{M}(\mathbb{N})$, the
characteristic function of $\mu$ is denoted by $\varphi_{\mu}$,
i.e. $\varphi_{\mu}(x)=\int e^{ixy}\: d\mu(y)$ for each $x\in\mathbb{R}$.

The following members of $\mathcal{M}(\mathbb{N})$ will appear later
on.

\subsubsection{The Poisson distribution}

For $0<t\in\mathbb{R}$ denote by $Pois(t)\in\mathcal{M}(\mathbb{N})$
the Poisson distribution with parameter $t$, which satisfies
\[
Pois(t)\{k\}=e^{-t}\frac{t^{k}}{k!}
\]
for each $k\in\mathbb{N}$.

\subsubsection{The Geometric distribution}

For $p\in[0,1)$ denote by $Geo(p)\in\mathcal{M}(\mathbb{N})$ the
geometric distribution with success parameter $p$, which satisfies
\[
Geo(p)\{k\}=(1-p)p^{k-1}
\]
for each $k\in\mathbb{N}_{+}$.

\subsubsection{The Compound Poisson distribution}

For $0<t\in\mathbb{R}$ and $\nu\in\mathcal{M}(\mathbb{N})$, denote
by $CP(t,\nu)\in\mathcal{M}(\mathbb{N})$ the compound Poisson distribution
with parameters $t$ and $\nu$, which satisfies
\[
CP(t,\nu)\{k\}=\sum_{j=1}^{\infty}Pois(t)\{j\}\cdot(\nu*)^{j}\{k\}
\]
for each $k\in\mathbb{N}$, where $(\nu*)^{j}$ is the $j$-fold convolution
of $\nu$. Let $W\sim Pois(t)$, and let $\eta_{1},\eta_{2},...$
be i.i.d random variables independent of $W$ with $\eta_{1}\sim\nu$,
then $\sum_{j=1}^{W}\eta_{j}\sim CP(t,\nu)$. Also, one checks that
\begin{equation}
\varphi_{CP(t,\nu)}(x)=\exp(t(\varphi_{\nu}(x)-1))\label{E8}
\end{equation}
for each $x\in\mathbb{R}$.

\subsubsection{The Pólya\textendash{}Aeppli distribution}

For $0<t\in\mathbb{R}$ and $p\in[0,1)$, denote by $PA(t,p)\in\mathcal{M}(\mathbb{N})$
the Pólya\textendash{}Aeppli distribution with parameters $t$ and
$p$, which satisfies
\[
PA(t,p)\{k\}=e^{-t}\sum_{j=1}^{k}\binom{k-1}{j-1}\frac{t^{j}}{j!}p^{k-j}(1-p)^{j}
\]
for each $k\in\mathbb{N}_{+}$, and $PA(t,p)\{0\}=e^{-t}$. One checks
that $PA(t,p)=CP(t,Geo(p))$, and so from (\ref{E8}),
\begin{equation}
\varphi_{PA(t,p)}(x)=\exp(t(\varphi_{Geo(p)}(x)-1))\label{E9}
\end{equation}
for each $x\in\mathbb{R}$. Observe that $PA(t,0)=Pois(t)$.

\section{\label{S3}Statement Of The Results}

\subsection{\label{SS5}Pointwise conditions for convergence and nonconvergence}

Throughout this article $t>0$ will be a fixed parameter. For each
$\omega\in\Omega_{\mathbb{P}}$ and $n\geq1$, let $N_{n}^{\omega}=[t(\mathbb{P}(A_{n}^{\omega}))^{-\ell}]$.
For each $r\in\mathbb{N}_{+}$ set
\[
\kappa(r)=lcm\{\frac{r}{gcd\{r,d_{j}\}}\::\:1\leq j\leq\ell\}
\]
where $lcm$ and $gcd$ denote the least common multiple and the greatest
common divisor, respectively. For an $n$-cylinder $A=[a_{0},...,a_{n-1}]$
with $\mathbb{P}(A)>0$, $r=\pi(A)$ and $R=[a_{0},...,a_{r-1}]$,
set
\[
\rho_{A}=\prod_{j=1}^{\ell}\mathbb{P}\{R^{(n+d_{j}\kappa(r))/r}\mid A\}.
\]
From (2.12) of Theorem 2.3 in \cite{key-1}, it follows that 
\begin{equation}
\sup\{\rho_{A}\::\: n\geq1,\: A\mbox{ is an \ensuremath{n}-cylinder, }\mathbb{P}(A)>0\}<1\label{E31}
\end{equation}
For a periodic point $\omega\in\Omega_{\mathbb{P}}$ with minimal
period $r\geq1$ (i.e. 
\[
r=\inf\{j\geq1\::\: T^{j}\omega=\omega\}
\]
and $r<\infty$) define $\beta_{\omega,n}=\mathbb{P}\{A_{n+r}^{\omega}\mid A_{n}^{\omega}\}$
for each $n\geq1$. The following simple lemma will be proven in Section
\ref{S4}.
\begin{lem}
\label{L1}Let $\omega\in\Omega_{\mathbb{P}}$ be a periodic point
with minimal period $r\geq1$. Assume that the limit $\beta_{\omega}=\underset{n\rightarrow\infty}{\lim}\beta_{\omega,n}$
exists. Then the limit $\rho_{\omega}=\underset{n\rightarrow\infty}{\lim}\rho_{A_{n}^{\omega}}$
exists, and is equal to $\beta_{\omega}^{a}$, where $a=\frac{\kappa(r)}{r}\sum_{i=1}^{\ell}d_{i}$.
Also, it holds that $\rho_{\omega}<1$.
\end{lem}

The pointwise conditions can now be stated.
\begin{thm}
\label{T2}Let $\omega\in\Omega_{\mathbb{P}}$ be a periodic point
with a minimal period $r\geq1$, then:

(a) If $\underset{n\rightarrow\infty}{\lim}\rho_{A_{n}^{\omega}}$
does not exists then $\{S_{N_{n}^{\omega}}^{A_{n}^{\omega}}\}_{n=1}^{\infty}$
does not converge in distribution.

(b) If the limit $\beta_{\omega}=\underset{n\rightarrow\infty}{\lim}\beta_{\omega,n}$
exists then $\mathcal{L}(S_{N_{n}^{\omega}}^{A_{n}^{\omega}})\overset{d}{\Longrightarrow}PA(t(1-\rho_{\omega}),\rho_{\omega})$
as $n\rightarrow\infty$.
\end{thm}

\begin{rem*}
Assertion (b) was first proven for the conventional setup in \cite{key-13}.
\end{rem*}

\begin{rem*}
Note that if $\rho_{\omega}=0$ then $\mathcal{L}(S_{N_{n}^{\omega}}^{A_{n}^{\omega}})\overset{d}{\Longrightarrow}Pois(t)$.
In Section \ref{SS3} an example in which this situation occurs will
be presented.
\end{rem*}

\begin{rem*}
Observe that in the conventional setup Theorem \ref{T2} says that
$\beta_{\omega}=\underset{n\rightarrow\infty}{\lim}\beta_{\omega,n}$
exists if and only if $\{S_{N_{n}^{\omega}}^{A_{n}^{\omega}}\}_{n=1}^{\infty}$
converges in distribution, in which case 
\[
\mathcal{L}(S_{N_{n}^{\omega}}^{A_{n}^{\omega}})\overset{d}{\Longrightarrow}PA(t(1-\beta_{\omega}),\beta_{\omega})\mbox{ as }n\rightarrow\infty.
\]

\end{rem*}

\subsection{\label{SS1}A global condition for convergence and applications}

\subsubsection{Continuous inverse Jacobian}

A condition on the dynamical system will now be stated that guarantees
the convergence in distribution of $\{S_{N_{n}^{\omega}}^{A_{n}^{\omega}}\}_{n=1}^{\infty}$,
for every periodic point $\omega\in\Omega_{\mathbb{P}}$.
\begin{thm}
\label{T3}Assume that the inverse Jacobian $J=\frac{d\mathbb{P}}{d\mathbb{P}\circ T}$
is a continuous function on $\Omega_{S}$ (i.e. it has a continuous
version), where $\Omega_{S}$ and its topology were defined in Section
\ref{SS2}. Let $\omega\in\Omega_{\mathbb{P}}$ be a periodic point
with minimal period $r\geq1$, then the limit $\underset{n\rightarrow\infty}{\lim}\beta_{\omega,n}$
exists and it is equal to $\prod_{j=0}^{r-1}J(T^{j}\omega)$.
\end{thm}

The next corollary follows immediately from Theorems \ref{T2} and
\ref{T3}.
\begin{cor}
\label{C1}Assume that the inverse Jacobian $J=\frac{d\mathbb{P}}{d\mathbb{P}\circ T}$
is a continuous function on $\Omega_{S}$. Let $\omega\in\Omega_{\mathbb{P}}$
be a periodic point with minimal period $r\geq1$, then
\[
\mathcal{L}(S_{N_{n}^{\omega}}^{A_{n}^{\omega}})\overset{d}{\Longrightarrow}PA(t(1-\rho_{\omega}),\rho_{\omega})\mbox{ as }n\rightarrow\infty
\]
where $\rho_{\omega}=\left(\prod_{j=0}^{r-1}J(T^{j}\omega)\right)^{a}$
and $a=\frac{\kappa(r)}{r}\sum_{i=1}^{\ell}d_{i}$.
\end{cor}

Two applications of Corollary \ref{C1} will now be stated.

\subsubsection{Gibbs measures}

In this section it will be assumed that the shift space is of finite
type, i.e. that $|\mathcal{A}|<\infty$. Let $\phi:\Omega_{S}\rightarrow\mathbb{R}$
be Hölder continuous with respect to the metric $d$, where $d$ was
defined in Section \ref{SS2}. From Theorem 1.4 in \cite{key-12}
it follows that there exist a unique $T$-invariant Borel probability
measure $\mathbb{P}_{\phi}$ on $\Omega_{S}$ for which one can find
constants $c_{1}>0$, $c_{2}>0$ and $P\in\mathbb{R}$ such that
\begin{equation}
c_{1}\leq\frac{\mathbb{P}_{\phi}[\omega_{0},...,\omega_{n-1}]}{\exp(-Pn+\sum_{j=0}^{n-1}\phi(T^{j}\omega))}\leq c_{2}\label{E12}
\end{equation}
for each $\omega\in\Omega_{S}$ and $n\geq0$. The measure $\mathbb{P}_{\phi}$
is called the Gibbs measure of $\phi$ and $P$ is called the pressure
of $\phi$. It is shown in \cite{key-12} that $\mathbb{P}_{\phi}$
is $\psi$-mixing. For each $E\in\mathcal{F}$ define $\mathbb{P}(E)=\mathbb{P}_{\phi}(E\cap\Omega_{S})$,
then $\mathbb{P}$ is as described in Section \ref{SS2}. Observe
that from (\ref{E12}) it follows that $\Omega_{\mathbb{P}}=\Omega_{S}$.
\begin{thm}
\label{T4}Let $J=\frac{d\mathbb{P}}{d\mathbb{P}\circ T}$ be the
inverse Jacobian of $\mathbb{P}$. Then there exist a continuous function
$h:\Omega_{S}\rightarrow\mathbb{R}$ with $h>0$, such that 
\[
J(\omega)=e^{\phi(\omega)-P}\frac{h(\omega)}{h(T\omega)}
\]
for $\mathbb{P}\circ T$-almost all $\omega\in\Omega_{S}$.
\end{thm}

From Theorem \ref{T4} and Corollary \ref{C1} it follows that:
\begin{cor}
\label{C2}Let $\omega\in\Omega_{S}$ be a periodic point with a minimal
period $r\geq1$, then
\[
\mathcal{L}(S_{N_{n}^{\omega}}^{A_{n}^{\omega}})\overset{d}{\Longrightarrow}PA(t(1-\rho_{\omega}),\rho_{\omega})\mbox{ as }n\rightarrow\infty
\]
where $\rho_{\omega}=\exp(a\cdot\sum_{j=0}^{r-1}(\phi(T^{j}\omega)-P))$
and $a$ is as in Corollary \ref{C1}.
\end{cor}

\begin{rem*}
For the conventional case Corollary \ref{C2} was proven in \cite{key-13}.
When $\mathbb{P}$ is a Bernoulli shift Corollary \ref{C2} was proven
in \cite{key-1} even for a countable alphabet.
\end{rem*}

\begin{rem*}
Corollary \ref{C2} actually still holds for some Gibbs measures over
a countable alphabet. In this case some additional assumptions should
be made in order to insure the $\psi$-mixing property.
\end{rem*}

\subsubsection{Expanding Markov Interval Maps}

Another family of systems, for which Corollary \ref{C1} can be applied,
will now be described. These systems are derived from certain expanding
Markov interval maps. In \cite{key-2} a family of dynamical systems
is defined, namely the family of $f$-expansions with their absolutely
continuous invariant measures, whose members satisfy all of the assumptions
listed below. An example of such a system is the Gauss map, $x\rightarrow\frac{1}{x}\:(\mbox{mod }1)$
for $x\in(0,1]$, equipped with the Gauss measure $\mu_{G}(\Gamma)=\frac{1}{\ln2}\int_{\Gamma}\frac{1}{1+x}\: dx$.
If we use the measure $\mu_{G}$ then the results from \cite{key-1}
could not be applied to numbers having periodic continued fraction
expansions (for instance the fractional parts of the golden ratio
$\frac{1+\sqrt{5}}{2}$ or $\sqrt{2}$), while Corollary \ref{C3}
below yields the Pólya\textendash{}Aeppli limiting distribution in
this case.

Let $I=[0,1]$ and let $m$ denote the Lebesgue measure on $I$. Let
$\{U_{a}\subset I\::\: a\in\mathcal{A}\}$ be a collection of disjoint
open intervals. For each $a\in\mathcal{A}$ set $I_{a}=\overline{U_{a}}$,
set $U=\bigcup_{a\in\mathcal{A}}U_{a}$ and assume that $m(I\setminus U)=0$.
Let $f:I\rightarrow I$ be such that: 

(i) For every $a\in\mathcal{A}$ the restriction of $f$ to $U_{a}$,
can be extended to a function $f_{a}\in C^{1}(I_{a})$ which is strictly
monotonic.

(ii) For every $a,b\in\mathcal{A}$, if $f(U_{a})\cap U_{b}\ne\emptyset$
then $U_{b}\subset f(U_{a})$.

(iii) For some integer $q\geq1$ there is a $\tau>1$, such that $|(f^{q})'(x)|\geq\tau$
for all $x\in\bigcap_{j=0}^{q-1}f^{-j}(U)$.

(iv) For every $a\in\mathcal{A}$ there are $b,c\in\mathcal{A}$ such
that $U_{b}\subset f(U_{a})$ and $U_{a}\subset f(U_{c})$.

Let $\mu$ be an $f$-invariant Borel probability measure on $I$
such that $\mu\ll m$. Set $p=\frac{d\mu}{dm}$, and assume that $p$
is continuous and strictly positive.

Set $N=\bigcup_{j=0}^{\infty}f^{-j}(I\setminus U)$ and $\tilde{I}=I\setminus N$.
The set $N$ is $\mu$-negligible, since
\[
\mu(N)\leq\sum_{j=0}^{\infty}\mu(f^{-j}(I\setminus U))=\sum_{j=0}^{\infty}\mu(I\setminus U)=\sum_{j=0}^{\infty}\underset{I\setminus U}{\int}\: p\: dm=0
\]
For each $j\geq0$, let $\xi_{j}:\tilde{I}\rightarrow\mathcal{A}$
be such that
\[
\xi_{j}(x)=a\mbox{ if and only if }f^{j}(x)\in U_{a}
\]
for each $x\in\tilde{I}$ and $a\in\mathcal{A}$. It is assumed that
the process $\{\xi_{j}\}_{j=0}^{\infty}$ is $\psi$-mixing (see (1.8)
in \cite{key-2}).

Define the matrix $S$, such that $S_{a,b}=1$ if and only if $U_{b}\subset f(U_{a})$,
for each $a,b\in\mathcal{A}$. We assumed that $T:\Omega_{S}\rightarrow\Omega_{S}$
is topologically mixing. It is not difficult to show (see Proposition
1.2 in \cite{key-10}) that there exist a Hölder continuous map $\Theta:\Omega_{S}\rightarrow I$
such that: 

(1) for every $x\in\tilde{I}$ there is a unique $\omega\in\Omega_{S}$
with $\Theta(\omega)=x$, and this $\omega$ satisfies $\omega_{j}=\xi_{j}(x)$
for each $j\geq1$.

(2) $\Theta(\omega)\in I_{\omega_{0}}$ for each $\omega\in\Omega_{S}$.

(3) if $\omega\in\Omega_{S}$ and $\Theta(\omega)\in\tilde{I}$ then
$\Theta(T\omega)=f(\Theta(\omega))$. 

For each $\omega\in\Omega_{S}$, $\Theta(\omega)$ is defined in \cite{key-10}
to be the unique element in $\bigcap_{n\geq1}\overline{\bigcap_{j=0}^{n-1}f^{-j}(U_{\omega_{j}})}$.

Let $\Phi:\tilde{I}\rightarrow\Omega$ be such that 
\[
\Phi(x)=(\xi_{0}(x),\xi_{1}(x),...)
\]
for each $x\in\tilde{I}$. Then $\Phi$ is the inverse of $\Theta$
when $\Theta$ is restricted to $\Theta^{-1}(\tilde{I})$. Let $\mathbb{P}=\mu\circ\Phi^{-1}$
then, since $\mu$ is $f$-invariant and $\{\xi_{j}\}_{j=0}^{\infty}$
is $\psi$-mixing, it follows that $\mathbb{P}$ is $T$-invariant
and $\psi$-mixing.
\begin{thm}
\label{T5}Let $J=\frac{d\mathbb{P}}{d\mathbb{P}\circ T}$ be the
inverse Jacobian of $\mathbb{P}$, then $J(\omega)=\frac{p(\Theta(\omega))}{f_{\omega_{0}}'(\Theta(\omega))p(\Theta(T\omega))}$
for $\mathbb{P}\circ T$-almost all $\omega\in\Omega_{S}$.
\end{thm}

From Theorem \ref{T5} and from Corollary \ref{C1} it follows that:
\begin{cor}
\label{C3}Let $\omega\in\Omega_{\mathbb{P}}$ be a periodic point
with minimal period $r\geq1$, then
\[
\mathcal{L}(S_{N_{n}^{\omega}}^{A_{n}^{\omega}})\overset{d}{\Longrightarrow}PA(t(1-\rho_{\omega}),\rho_{\omega})\mbox{ as }n\rightarrow\infty
\]
where $\rho_{\omega}=\left(\prod_{j=0}^{r-1}f_{\omega_{j}}'(\Theta(T^{j}\omega))\right)^{-a}$,
and $a$ is as in Corollary \ref{C1}.
\end{cor}

\subsection{Analysis of partial limits}

As we can see from the examples constructed in Section \ref{SS4},
for a periodic point $\omega$ the sequence $\{S_{N_{n}^{\omega}}^{A_{n}^{\omega}}\}$
does not necessarily has a limit in distribution. But its partial
limits are always Compound Poisson. 
\begin{thm}
\label{T6} Let $\omega\in\Omega_{\mathbb{P}}$ be a periodic point,
and for each $n\geq1$ let $\mu_{n}\in\mathcal{M}(\mathbb{N})$ be
the distribution of $S_{N_{n}^{\omega}}^{A_{n}^{\omega}}$. Then:

(a) The sequence $\{\mu_{n}\}$ is tight.

(b) Let $\mu$ be a probability distribution on $\mathbb{R}$ such
that $\mu_{n_{k}}\overset{d}{\Longrightarrow}\mu$ as $k\rightarrow\infty$,
for some sub-sequence $\{\mu_{n_{k}}\}$. Then $\mu=CP(\tau,\theta)$,
for some $0<\tau\leq t$ and $\theta\in\mathcal{M}(\mathbb{N})$.
\end{thm}

\subsection{\label{SS4}A family of nonconvergence examples}

In \cite{key-1} an example of a $\psi$-mixing system was built in
which there exists a periodic point $\omega$ such that $\{S_{N_{n}^{\omega}}^{A_{n}^{\omega}}\}$
does not converge in distribution. This system was derived from a
Bernoulli shift over the alphabet $\{0,1\}$ by using the group structure
of $\{0,1\}$. The next result extends this construction to every
finite abelian group, producing a family of nonconvergence examples.

Let $(G,+)$ be a finite abelian group, and let $(p_{g})_{g\in G}\subset(0,1)$
be such that $\sum_{g\in G}p_{g}=1$. Assume $\Omega=G^{\mathbb{N}}$
and let $\mathbb{P}$ be the measure on $(\Omega,\mathcal{F})$ that
satisfies
\[
\mathbb{P}[g_{0},...,g_{n-1}]=\prod_{j=0}^{n-1}p_{g_{j}}
\]
for each $g_{0},...,g_{n-1}\in G$. Then the coordinate projections
from $\Omega$ onto $G$ are i.i.d. random elements.

Let $N\geq2$ be an integer, let $\Phi:\Omega\rightarrow\Omega$ be
such that
\[
(\Phi\omega)_{j}=\omega_{j}+...+\omega_{j+N-1}
\]
for each $\omega\in\Omega$ and $j\in\mathbb{N}$, and set $\mathbb{P}_{0}=\mathbb{P}\circ\Phi^{-1}$.
\begin{thm}
\label{T7}(a) $\mathbb{P}_{0}$ is $T$-invariant.

(b) \noun{$\mathbb{P}_{0}$ }is $\psi$-mixing.

(c) Assume there exist $h\in G$ with $p_{h}>p_{g}$ for all $g\in G\setminus\{h\}$.
Then there exists $s\in G$ such that for $\omega\in\Omega$ with
$\omega_{j}=s$ for each $j\geq0$, the limit $\underset{n}{\lim}\:\mathbb{P}_{0}\{A_{n+1}^{\omega}\mid A_{n}^{\omega}\}$
does not exist.
\end{thm}

From Theorem \ref{T7} we obtain:
\begin{cor}
(a) In the conventional setup, it follows from Assertion (a) of Theorem
\ref{T2} that $\{S_{N_{n}^{\omega}}^{A_{n}^{\omega}}\}$ does not
converge in distribution (where $\omega$ is as in Assertion (c) of
Theorem \ref{T7}).

(b) It follows from Theorem \ref{T3} that the inverse Jacobian $\frac{d\mathbb{P}_{0}}{d\mathbb{P}_{0}\circ T}$
does not have a continuous version.

(c) It follows from Corollary \ref{C2} that $\mathbb{P}_{0}$ is
not a Gibbs measure corresponding to a Hölder continuous function.
\end{cor}

Hence our construction provides a large class of $\psi$-mixing non
Gibbs measures which seems to be new.

\subsection{\label{SS3}Example of a pure Poisson limit distribution}

An example of a system will now be constructed in which there exists
a periodic point $\omega$ such that the limit distribution of $\{S_{N_{n}^{\omega}}^{A_{n}^{\omega}}\}$
is pure Poisson.

Set $\mathcal{A}=\mathbb{N}_{+}$ and $\Omega=\mathcal{A}^{\mathbb{N}}$,
and let $\mathbb{P}$ be the measure on $(\Omega,\mathcal{F})$ such
that
\[
\mathbb{P}[a_{0},...,a_{n-1}]=\prod_{j=0}^{n-1}2^{-a_{j}}
\]
for each $a_{0},...,a_{n-1}\in\mathcal{A}$. Then the coordinate projections
from $\Omega$ onto $\mathcal{A}$ are i.i.d. random variables.

Let $\Omega_{0}=\{0,1\}^{\mathbb{N}}$, let $\mathcal{F}_{0}$ be
the $\sigma$-algebra generated by cylinders, and let $T_{0}:\Omega_{0}\rightarrow\Omega_{0}$
be the left-shift operator. For each $a_{1},a_{2}\in\mathcal{A}$
set
\[
\theta(a_{1},a_{2})=\begin{cases}
1 & ,\mbox{ if }a_{2}=a_{1}+1\\
0 & ,\mbox{ otherwise}
\end{cases}
\]
Let $\Theta:\Omega\rightarrow\Omega_{0}$ be such that $(\Theta\omega)_{j}=\theta(\omega_{j},\omega_{j+1})$
for each $\omega\in\Omega$ and $j\geq0$, and let $\mathbb{P}_{0}=\mathbb{P}\circ\Theta^{-1}$.
\begin{thm}
\label{T8}(a) $\mathbb{P}_{0}$ is $T_{0}$-invariant.

(b) \noun{$\mathbb{P}_{0}$ }is $\psi$-mixing.

(c) Let $\omega\in\Omega_{0}$ be such that $\omega_{j}=1$ for each
$j\geq1$, then $\underset{n}{\lim}\:\mathbb{P}_{0}\{A_{n+1}^{\omega}\mid A_{n}^{\omega}\}=0$.
\end{thm}

From Theorem \ref{T8} we obtain:
\begin{cor}
(a) From Assertion (b) of Theorem \ref{T2} it follows that $\mathcal{L}(S_{N_{n}^{\omega}}^{A_{n}^{\omega}})\overset{d}{\Longrightarrow}Pois(t)$
as $n\rightarrow\infty$, where $\omega$ is as in Assertion (c) of
Theorem \ref{T8}.

(b) From Corollary \ref{C2} it follows that $\mathbb{P}_{0}$ is
not a Gibbs measure corresponding to a Hölder continuous function.
\end{cor}

\section{\label{S4}Proof of Theorem \ref{T2}}

\emph{Proof of Lemma \ref{L1}:} It holds that
\begin{multline*}
\underset{n\rightarrow\infty}{\lim}\rho_{A_{n}^{\omega}}=\underset{n\rightarrow\infty}{\lim}\prod_{i=1}^{\ell}\mathbb{P}\{R^{(n+d_{i}\kappa(r))/r}\mid R^{n/r}\}=\\
=\underset{n\rightarrow\infty}{\lim}\prod_{i=1}^{\ell}\prod_{j=1}^{d_{i}\kappa(r)/r}\frac{\mathbb{P}(R^{(n+jr)/r})}{\mathbb{P}(R^{(n+(j-1)r)/r})}=\underset{n\rightarrow\infty}{\lim}\prod_{i=1}^{\ell}\prod_{j=1}^{d_{i}\kappa(r)/r}\beta_{\omega,n+(j-1)r}=\beta_{\omega}^{a}
\end{multline*}
and from (\ref{E31}) we see that $\underset{n\rightarrow\infty}{\lim}\rho_{A_{n}^{\omega}}<1$.
$\square$

\emph{Proof of Theorem \ref{T2}: }For each $n\geq r$ let $W_{n}\sim Pois(t(1-\rho_{A_{n}^{\omega}}))$.
From Theorem 2.3 in \cite{key-1} it follows that for each $n\geq1$
there exists a sequence of i.i.d. random variables $\eta_{n,1},\eta_{n,2},...$
independent of $W_{n}$, such that $\mathbb{P}\{\eta_{n,1}\in\{1,...,[\frac{n}{r}]\}\}=1$
and such that for $Y_{n}=\sum_{j=1}^{W_{n}}\eta_{n,j}$
\begin{equation}
\underset{n\rightarrow\infty}{\lim}d_{TV}(\mathcal{L}(S_{N_{n}^{\omega}}^{A_{n}^{\omega}}),\mathcal{L}(Y_{n}))=0\label{E1}
\end{equation}

\emph{(a)} Assume that $\underset{n\rightarrow\infty}{\lim}\rho_{A_{n}^{\omega}}$
does not exist, then there exists $\epsilon>0$ such that for each
$M\geq1$, there exist $n>m>M$ with $|\rho_{A_{n}^{\omega}}-\rho_{A_{m}^{\omega}}|>\epsilon$.
Let $M\geq1$ be such that $d_{TV}(\mathcal{L}(S_{N_{n}^{\omega}}^{A_{n}^{\omega}}),\mathcal{L}(Y_{n}))<\frac{e^{-t}t\epsilon}{3}$
for each $n>M$, then for each $n>m>M$ with $|\rho_{A_{n}^{\omega}}-\rho_{A_{m}^{\omega}}|>\epsilon$,
\begin{multline*}
|\mathbb{P}\{S_{N_{n}^{\omega}}^{A_{n}^{\omega}}=0\}-\mathbb{P}\{S_{N_{m}^{\omega}}^{A_{m}^{\omega}}=0\}|\geq\\
\geq|\mathbb{P}\{Y_{n}=0\}-\mathbb{P}\{Y_{m}=0\}|-d_{TV}(\mathcal{L}(S_{N_{n}^{\omega}}^{A_{n}^{\omega}}),\mathcal{L}(Y_{n}))-d_{TV}(\mathcal{L}(S_{N_{m}^{\omega}}^{A_{m}^{\omega}}),\mathcal{L}(Y_{m}))>\\
>|\exp(-t(1-\rho_{A_{n}^{\omega}}))-\exp(-t(1-\rho_{A_{m}^{\omega}}))|-\frac{e^{-t}t\epsilon}{3}-\frac{e^{-t}t\epsilon}{3}\geq\\
\geq e^{-t}t|\rho_{A_{n}^{\omega}}-\rho_{A_{m}^{\omega}}|-2\frac{e^{-t}t\epsilon}{3}>\frac{e^{-t}t\epsilon}{3}
\end{multline*}
which shows that $\{S_{N_{n}^{\omega}}^{A_{n}^{\omega}}\}_{n=1}^{\infty}$
does not converge in distribution.

\emph{(b)} Assume that $\beta_{\omega}=\underset{n\rightarrow\infty}{\lim}\beta_{\omega,n}$
exists, then from Lemma \ref{L1} it follows that the limit $\rho_{\omega}=\underset{n\rightarrow\infty}{\lim}\rho_{A_{n}^{\omega}}$
exists, it is equal to $\beta_{\omega}^{a}$, and is strictly less
than $1$. Because of (\ref{E1}) it is enough to prove that $\mathcal{L}(Y_{n})\overset{d}{\Longrightarrow}PA(t(1-\rho_{\omega}),\rho_{\omega})$
as $n\rightarrow\infty$.

The main part of the proof will be showing that $\mathcal{L}(\eta_{n,1})\overset{d}{\Longrightarrow}Geo(\rho_{\omega})$
as $n\rightarrow\infty$. To do this fix some integer $b\geq1$. It
is enough to show that $\mathbb{P}\{\eta_{n,1}=b\}\overset{n\rightarrow\infty}{\longrightarrow}Geo(\rho_{\omega})\{b\}$. 

Since $\rho_{\omega}<1$ it follows that 
\[
1-2\beta_{\omega}^{a}+\beta_{\omega}^{2a}=1-2\rho_{\omega}+\rho_{\omega}^{2}>0.
\]
From this and from $\psi_{n}\overset{n\rightarrow\infty}{\longrightarrow}0$
we see that there exist $M\geq1$ and $\epsilon>0$ such that 
\begin{equation}
1-2\prod_{i=1}^{\ell}\prod_{j=1}^{d_{i}\kappa(r)/r}\beta_{\omega,n+(j-1)r}+\prod_{i=1}^{\ell}\prod_{j=1}^{2d_{i}\kappa(r)/r}\beta_{\omega,n+(j-1)r}-4\cdot2^{\ell}\psi_{n}>\epsilon\label{E3}
\end{equation}
for all $n>M$. By choosing $M$ large enough and using (\ref{E2})
we may assume that $N_{n}^{\omega}>15d_{\ell}rn$, $n>d_{\ell}\kappa(r)b$
and $\psi_{n}<(3/2)^{1/(\ell+1)}-1$ for each $n\geq M$.$\newline$$\newline$Let
$n\geq M$. We will use the notation from \cite{key-1} appearing
in the statement and the proof of Theorem 2.3 there, with $A=A_{n}^{\omega}$.
Namely, we set: $R=[\omega_{0},...,\omega_{r-1}]$, $N=N_{n}^{\omega}$,
$K=5d_{\ell}rn$, $S_{N}=S_{N}^{A}$, $X_{k}=X_{k}^{A}$ for each
$k\in\mathbb{N}_{+}$, $\hat{X}_{\alpha}=1_{\{K<\alpha\leq N\}}\cdot X_{\alpha}$
for each $\alpha\in\mathbb{Z}$, $n_{0}=[\frac{n}{r}]$, $\kappa=\kappa(r)$,
$\rho=\rho_{A}$, $I_{0}=\{K+1,...,N\}\times\{1,...,n_{0}\}$,
\[
X_{\alpha,j}=(1-\hat{X}_{\alpha-\kappa})(1-\hat{X}_{\alpha+j\kappa})\prod_{k=0}^{j-1}\hat{X}_{\alpha+k\kappa}\mbox{ for each }\alpha,j\in\mathbb{N}_{+},
\]
$\lambda_{\alpha,j}=E[X_{\alpha,j}]$ for each $\alpha,j\in\mathbb{N}_{+}$,
$\lambda=\sum_{(\alpha,j)\in I_{0}}\lambda_{\alpha,j}$, $\lambda_{j}=\lambda^{-1}\sum_{\alpha=K+1}^{N}\lambda_{\alpha,j}$
for each $1\leq j\leq n_{0}$, and $s=t(1-\rho)$. We recall that
in the proof of Theorem 2.3 in \cite{key-1}, the i.i.d. random variables
$\eta_{n,1},\eta_{n,2},...$ were chosen so that $\mathbb{P}\{\eta_{n,1}=j\}=\lambda_{j}$
for each $1\leq j\leq n_{0}$. First we shall need to bound $\lambda$
from below.\emph{$\newline$$\newline$Claim:} It holds that $\lambda>\frac{t\epsilon}{2}$.\emph{$\newline$$\newline$Proof:
}Let $K<\alpha\leq N$, then from Lemma 3.2 in \cite{key-1} it follows
that
\begin{multline*}
\lambda_{\alpha,1}=E[X_{\alpha,1}]\geq E[(1-X_{\alpha-\kappa})(1-X_{\alpha+\kappa})X_{\alpha}]=\\
=E[X_{\alpha}]-E[X_{\alpha-\kappa}X_{\alpha}]-E[X_{\alpha+\kappa}X_{\alpha}]+E[X_{\alpha-\kappa}X_{\alpha+\kappa}X_{\alpha}]=\\
=\mathbb{P}(\bigcap_{i=1}^{\ell}T^{-d_{i}\alpha}(A))-\mathbb{P}(\bigcap_{i=1}^{\ell}T^{-d_{i}(\alpha-\kappa)}(R^{(n+d_{i}\kappa)/r}))-\\
-\mathbb{P}(\bigcap_{i=1}^{\ell}T^{-d_{i}\alpha}(R^{(n+d_{i}\kappa)/r}))+\mathbb{P}(\bigcap_{i=1}^{\ell}T^{-d_{i}(\alpha-\kappa)}(R^{(n+2d_{i}\kappa)/r}))\geq\\
\geq(\mathbb{P}(A))^{\ell}(1-2^{\ell}\psi_{n})-2(\prod_{i=1}^{\ell}\mathbb{P}(R^{(n+d_{i}\kappa)/r}))(1+2^{\ell}\psi_{n})+(1-2^{\ell}\psi_{n})(\prod_{i=1}^{\ell}\mathbb{P}(R^{(n+2d_{i}\kappa)/r}))\geq\\
\geq(\mathbb{P}(A))^{\ell}-2\prod_{i=1}^{\ell}\mathbb{P}(R^{(n+d_{i}\kappa)/r})+\prod_{i=1}^{\ell}\mathbb{P}(R^{(n+2d_{i}\kappa)/r})-4\cdot2^{\ell}\psi_{n}(\mathbb{P}(A))^{\ell}.
\end{multline*}
Hence, from (\ref{E3}),
\begin{multline*}
(\mathbb{P}(A))^{-\ell}\lambda_{\alpha,1}\geq1-2\prod_{i=1}^{\ell}\frac{\mathbb{P}(R^{(n+d_{i}\kappa)/r})}{\mathbb{P}(A)}+\prod_{i=1}^{\ell}\frac{\mathbb{P}(R^{(n+2d_{i}\kappa)/r})}{\mathbb{P}(A)}-4\cdot2^{\ell}\psi_{n}=\\
=1-2\prod_{i=1}^{\ell}\prod_{j=1}^{d_{i}\kappa/r}\frac{\mathbb{P}(R^{(n+jr)/r})}{\mathbb{P}(R^{(n+(j-1)r)/r})}+\prod_{i=1}^{\ell}\prod_{j=1}^{2d_{i}\kappa/r}\frac{\mathbb{P}(R^{(n+jr)/r})}{\mathbb{P}(R^{(n+(j-1)r)/r})}-4\cdot2^{\ell}\psi_{n}=\\
=1-2\prod_{i=1}^{\ell}\prod_{j=1}^{d_{i}\kappa/r}\beta_{\omega,n+(j-1)r}+\prod_{i=1}^{\ell}\prod_{j=1}^{2d_{i}\kappa/r}\beta_{\omega,n+(j-1)r}-4\cdot2^{\ell}\psi_{n}>\epsilon,
\end{multline*}
and so
\[
\lambda\geq\sum_{\alpha=K+1}^{N}\lambda_{\alpha,1}\geq(N-K)(\mathbb{P}(A))^{\ell}\epsilon>\frac{t\epsilon}{2}.
\]

We resume the proof of Assertion (b) of Theorem \ref{T2}. From $\lambda^{-1}\leq\frac{2}{t\epsilon}$
it follows that
\begin{multline*}
|\mathbb{P}\{\eta_{n,1}=b\}-Geo(\rho_{\omega})\{b\}|=|\lambda_{b}-(1-\rho_{\omega})\rho_{\omega}^{b-1}|=\\
=\lambda^{-1}|\sum_{\alpha=K+1}^{N}\lambda_{\alpha,b}-\lambda(1-\rho_{\omega})\rho_{\omega}^{b-1}|\leq\frac{2}{t\epsilon}|\sum_{\alpha=K+1}^{N}\lambda_{\alpha,b}-s(1-\rho_{\omega})\rho_{\omega}^{b-1}|+\frac{2}{t\epsilon}|\lambda-s|.
\end{multline*}
From the inequality (5.5) in the proof of Theorem 2.3 in \cite{key-1},
it follows that $|\lambda-s|$ tends to $0$ as $n\rightarrow\infty$.
Hence in order to prove that $\mathcal{L}(\eta_{n,1})\overset{d}{\Longrightarrow}Geo(\rho_{\omega})$,
it is enough to show that
\[
\underset{n\rightarrow\infty}{\lim}|\sum_{\alpha=K+1}^{N}\lambda_{\alpha,b}-s(1-\rho_{\omega})\rho_{\omega}^{b-1}|=0.
\]
Also, since
\[
s=t(1-\rho_{A_{n}^{\omega}})\overset{n\rightarrow\infty}{\longrightarrow}t(1-\rho_{\omega})
\]
it is enough to show that
\[
\underset{n\rightarrow\infty}{\lim}|\sum_{\alpha=K+1}^{N}\lambda_{\alpha,b}-t(1-\rho_{\omega})^{2}\rho_{\omega}^{b-1}|=0.
\]
We also have
\begin{multline*}
|\sum_{\alpha=K+1}^{N}\lambda_{\alpha,b}-t(1-\rho_{\omega})^{2}\rho_{\omega}^{b-1}|\leq\\
\leq(\mathbb{P}(A))^{\ell}+|\sum_{\alpha=K+1}^{N}\lambda_{\alpha,b}-N(\mathbb{P}(A))^{\ell}(1-\rho_{\omega})^{2}\rho_{\omega}^{b-1}|\leq\\
\leq5K\mathbb{P}(A)+\sum_{\alpha=K+\kappa+1}^{N-b\kappa}|\lambda_{\alpha,b}-(\mathbb{P}(A))^{\ell}(1-\rho_{\omega})^{2}\rho_{\omega}^{b-1}|,
\end{multline*}
and from (\ref{E2}), 
\[
K\mathbb{P}(A)\leq5d_{\ell}rn\cdot e^{-\Gamma n}\overset{n\rightarrow\infty}{\longrightarrow}0.
\]
Hence, it is enough to show that
\[
\underset{n\rightarrow\infty}{\lim}\sum_{\alpha=K+\kappa+1}^{N-b\kappa}|\lambda_{\alpha,b}-(\mathbb{P}(A))^{\ell}(1-\rho_{\omega})^{2}\rho_{\omega}^{b-1}|=0.
\]
This will be established once we prove the following:$\newline$$\newline$\emph{Claim:}
It holds that
\begin{equation}
\underset{n\rightarrow\infty}{\lim}(\mathbb{P}(A))^{-\ell}\sup\{|\lambda_{\alpha,b}-(\mathbb{P}(A))^{\ell}(1-\rho_{\omega})^{2}\rho_{\omega}^{b-1}|\::\: K+\kappa<\alpha\leq N-b\kappa\}=0.\label{E7}
\end{equation}
\emph{Proof:} Let $K+\kappa<\alpha\leq N-b\kappa$, then
\begin{multline}
|\lambda_{\alpha,b}-(\mathbb{P}(A))^{\ell}(1-\rho_{\omega})^{2}\rho_{\omega}^{b-1}|=\\
=|E[(1-X_{\alpha-\kappa})(1-X_{\alpha+b\kappa})\prod_{k=0}^{b-1}X_{\alpha+k\kappa}]-(\mathbb{P}(A))^{\ell}(1-\rho_{\omega})^{2}\rho_{\omega}^{b-1}|=\\
\leq|E[\prod_{k=0}^{b-1}X_{\alpha+k\kappa}]-(\mathbb{P}(A))^{\ell}\rho_{\omega}^{b-1}|+|E[\prod_{k=-1}^{b-1}X_{\alpha+k\kappa}]-(\mathbb{P}(A))^{\ell}\rho_{\omega}^{b}|+\\
+|E[\prod_{k=0}^{b}X_{\alpha+k\kappa}]-(\mathbb{P}(A))^{\ell}\rho_{\omega}^{b}|+|E[\prod_{k=-1}^{b}X_{\alpha+k\kappa}]-(\mathbb{P}(A))^{\ell}\rho_{\omega}^{b+1}|=\\
=\delta(\alpha,0,b-1)+\delta(\alpha,-1,b-1)+\delta(\alpha,0,b)+\delta(\alpha,-1,b)\label{E4}
\end{multline}
where
\[
\delta(\alpha,q,p)=|E[\prod_{k=q}^{p}X_{\alpha+k\kappa}]-(\mathbb{P}(A))^{\ell}\rho_{\omega}^{p-q}|
\]
for each $K+\kappa<\alpha\leq N-b\kappa$, $q\in\{-1,0\}$ and $p\in\{b-1,b\}$.
Fix such $\alpha$, $q$ and $p$, then 
\begin{multline}
\delta(\alpha,q,p)\leq|E[\prod_{k=q}^{p}X_{\alpha+k\kappa}]-\prod_{i=1}^{\ell}\mathbb{P}(R^{(n+(p-q)d_{i}\kappa)/r})|+\\
+|\prod_{i=1}^{\ell}\mathbb{P}(R^{(n+(p-q)d_{i}\kappa)/r})-(\mathbb{P}(A))^{\ell}\prod_{i=1}^{\ell}\prod_{k=1}^{p-q}\mathbb{P}\{R^{(n+kd_{i}\kappa)/r}\mid R^{(n+(k-1)d_{i}\kappa)/r}\}|+\\
+|(\mathbb{P}(A))^{\ell}\prod_{i=1}^{\ell}\prod_{k=1}^{p-q}\mathbb{P}\{R^{(n+kd_{i}\kappa)/r}\mid R^{(n+(k-1)d_{i}\kappa)/r}\}-(\mathbb{P}(A))^{\ell}\rho_{\omega}^{p-q}|=\Lambda_{1}+\Lambda_{2}+\Lambda_{3}\label{E5}
\end{multline}
where $\Lambda_{1}$, $\Lambda_{2}$ and $\Lambda_{3}$ denote the
first, second and third terms respectively, and if $p=q$ then the
multiplication from $1$ to $p-q$ equals $1$ by definition. From
Lemma 3.2 in \cite{key-1} it follows that
\begin{multline*}
\Lambda_{1}=|\mathbb{P}(\bigcap_{i=1}^{\ell}T^{-(\alpha+q\kappa)d_{i}}(R^{(n+(p-q)d_{i}\kappa)/r}))-\prod_{i=1}^{\ell}\mathbb{P}(R^{(n+(p-q)d_{i}\kappa)/r})|\leq\\
\leq((1+\psi_{n})^{\ell}-1)\prod_{i=1}^{\ell}\mathbb{P}(R^{(n+(p-q)d_{i}\kappa)/r})\leq2^{\ell}\psi_{n}(\mathbb{P}(A))^{\ell}.
\end{multline*}
The term $\Lambda_{2}$ vanishes and
\begin{multline*}
\Lambda_{3}=(\mathbb{P}(A))^{\ell}|\prod_{i=1}^{\ell}\prod_{k=1}^{p-q}\;\prod_{j=1+(k-1)d_{i}\kappa/r}^{kd_{i}\kappa/r}\frac{\mathbb{P}(R^{(n+jr)/r})}{\mathbb{P}(R^{(n+(j-1)r)/r})}-\rho_{\omega}^{p-q}|=\\
=(\mathbb{P}(A))^{\ell}|\prod_{i=1}^{\ell}\prod_{k=1}^{p-q}\;\prod_{j=1+(k-1)d_{i}\kappa/r}^{kd_{i}\kappa/r}\beta_{\omega,n+(j-1)r}-\rho_{\omega}^{p-q}|.
\end{multline*}
Hence from (\ref{E4}) and (\ref{E5}) it follows that
\begin{multline}
(\mathbb{P}(A))^{-\ell}\sup\{|\lambda_{\alpha,b}-(\mathbb{P}(A))^{\ell}(1-\rho_{\omega})^{2}\rho_{\omega}^{b-1}|\::\: K+\kappa<\alpha\leq N-b\kappa\}\leq\\
\leq4(\mathbb{P}(A))^{-\ell}\sup\{\delta(\alpha,q,p)\::\: K+\kappa<\alpha\leq N-b\kappa,\: q\in\{-1,0\},\: p\in\{b-1,b\}\}\leq\\
\leq4\cdot2^{\ell}\psi_{n}+4|\prod_{i=1}^{\ell}\prod_{k=1}^{p-q}\;\prod_{j=1+(k-1)d_{i}\kappa/r}^{kd_{i}\kappa/r}\beta_{\omega,n+(j-1)r}-\rho_{\omega}^{p-q}|.\label{E6}
\end{multline}
From $\beta_{\omega}=\underset{n\rightarrow\infty}{\lim}\beta_{\omega,n}$
it follows that
\begin{multline*}
\underset{n\rightarrow\infty}{\lim}\prod_{i=1}^{\ell}\prod_{k=1}^{p-q}\;\prod_{j=1+(k-1)d_{i}\kappa/r}^{kd_{i}\kappa/r}\beta_{\omega,n+(j-1)r}=\\
=\prod_{i=1}^{\ell}\prod_{k=1}^{p-q}\;\prod_{j=1+(k-1)d_{i}\kappa/r}^{kd_{i}\kappa/r}\beta_{\omega}=\prod_{k=1}^{p-q}\prod_{i=1}^{\ell}\beta^{(d_{i}\kappa)/r}=\rho_{\omega}^{p-q}
\end{multline*}
so from $\psi_{n}\overset{n\rightarrow\text{\ensuremath{\infty}}}{\longrightarrow}0$
and from (\ref{E6}) we obtain (\ref{E7}).$\newline$$\newline$We
have thus proven that $\mathcal{L}(\eta_{n,1})\overset{d}{\Longrightarrow}Geo(\rho_{\omega})$
as $n\rightarrow\infty$, and we can now finally show that $\mathcal{L}(Y_{n})\overset{d}{\Longrightarrow}PA(t(1-\rho_{\omega}),\rho_{\omega})$
as $n\rightarrow\infty$. Let $x\in\mathbb{R}$, then from Theorem
26.3 in \cite{key-9} it follows that
\[
\underset{n\rightarrow\infty}{\lim}\varphi_{\eta_{n,1}}(x)=\varphi_{Geo(\rho_{\omega})}(x),
\]
where recall that we denote by $\varphi_{\xi}$ the characteristic
function of a random variable or measure $\xi$. Hence, by (\ref{E8})
and (\ref{E9}), 
\begin{multline*}
\underset{n\rightarrow\infty}{\lim}\varphi_{\mathcal{L}(Y_{n})}(x)=\underset{n\rightarrow\infty}{\lim}\varphi_{CP(t(1-\rho_{A_{n}^{\omega}}),\mathcal{L}(\eta_{n,1}))}(x)=\\
=\underset{n\rightarrow\infty}{\lim}\exp(t(1-\rho_{A_{n}^{\omega}})(\varphi_{\eta_{n,1}}(x)-1))=\\
=\exp(t(1-\rho_{\omega})(\varphi_{Geo(\rho_{\omega})}(x)-1))=\varphi_{PA(t(1-\rho_{\omega}),\rho_{\omega})}(x).
\end{multline*}
Now another application of Theorem 26.3 gives $\mathcal{L}(Y_{n})\overset{d}{\Longrightarrow}PA(t(1-\rho_{\omega}),\rho_{\omega})$
as $n\rightarrow\infty$, which completes the proof of the theorem.
$\square$

\section{\label{S5}Proof of the results of Section \ref{SS1}}

In this section we will write $T$ for the restriction of $T$ to
$\Omega_{S}$ and for each $w\in\mathcal{A}_{S}^{*}$ we will write
$[w]$ for $[w]\cap\Omega_{S}$ .

\emph{Proof of Theorem \ref{T3}:} Set $w=\omega_{0}\cdot...\cdot\omega_{r-1}$
and for each $n\geq1$ and $0\leq j<r$ set 
\[
E_{n,j}=[\omega_{r-1-j}\cdot...\cdot\omega_{r-1}\cdot w^{n/r}]
\]
then
\[
T(E_{n,j})=\begin{cases}
[w^{n/r}] & \mbox{, if }j=0\\
{}[\omega_{r-j}\cdot...\cdot\omega_{r-1}\cdot w^{n/r}] & \mbox{, if }0<j<r
\end{cases}.
\]
Hence, given $n\geq1$,
\begin{equation}
\beta_{\omega,n}=\frac{\mathbb{P}([w^{(n+r)/r}])}{\mathbb{P}([w^{n/r}])}=\prod_{j=0}^{r-1}\frac{\mathbb{P}(E_{n,j})}{\mathbb{P}(T(E_{n,j}))}=\prod_{j=0}^{r-1}\frac{1}{\mathbb{P}\circ T(E_{n,j})}\int_{E_{n,j}}J\: d\mathbb{P}\circ T.\label{E10}
\end{equation}

Let $0\leq j<r$ and $\epsilon>0$. Since $J$ is continuous, since
$T^{r-1-j}(\omega)\in E_{n,j}$ for all $n\geq1$, and since
\[
\sup\{d(\gamma,\eta)\::\:\gamma,\eta\in E_{n,j}\}\overset{n\rightarrow\infty}{\longrightarrow}0
\]
(where $d$ is the metric defined in section \ref{SS2}), it follows
that there exist $N\geq1$ with 
\[
|J(T^{r-1-j}(\omega))-J(\gamma)|<\epsilon
\]
for each $n\geq N$ and $\gamma\in E_{n,j}$. From this it follows
that for each $n\geq N$ 
\begin{multline*}
|\frac{1}{\mathbb{P}\circ T(E_{n,j})}\int_{E_{n,j}}J\: d\mathbb{P}\circ T-J(T^{r-1-j}(\omega))|\leq\\
\leq\frac{1}{\mathbb{P}\circ T(E_{n,j})}\int_{E_{n,j}}|J(\gamma)-J(T^{r-1-j}(\omega))|d\mathbb{P}\circ T(\gamma)\leq\epsilon,
\end{multline*}
and so
\begin{equation}
\underset{n}{\lim}\:\frac{1}{\mathbb{P}\circ T(E_{n,j})}\int_{E_{n,j}}J\: d\mathbb{P}\circ T=J(T^{r-1-j}(\omega)).\label{E11}
\end{equation}
The theorem now follows from (\ref{E10}) and (\ref{E11}). $\square$

\emph{Proof of Theorem \ref{T4}:} Let $C(\Omega_{S})$ be the Banach
space of all continuous functions from $\Omega_{S}$ to $\mathbb{R}$.
Let $L_{\phi}:C(\Omega_{S})\rightarrow C(\Omega_{S})$ be the Ruelle
(transfer) operator associated with $\phi$, i.e.
\[
(L_{\phi}f)(x)=\sum_{y\in T^{-1}\{x\}}e^{\phi(y)}f(y)
\]
for each $f\in C(\Omega_{S})$ and $x\in\Omega_{S}$. From the construction
of the Gibbs measure $\mathbb{P}_{\phi}$ carried out in sections
B and C of Chapter 1 in \cite{key-12}, it follows that there exist
$h\in C(\Omega_{S})$ with $h>0$ and a probability measure $\nu$
on $(\Omega_{S},\mathcal{F}_{S})$ such that $L_{\phi}h=e^{P}h$,
$L_{\phi}^{*}\nu=e^{P}\nu$ and $d\mathbb{P}_{\phi}=h\: d\nu$. Here
$L_{\phi}^{*}$ is the adjoint operator of $L_{\phi}$, which satisfies
\[
\int f\: d(L_{\phi}^{*}\mu)=\int L_{\phi}f\: d\mu
\]
for each finite measure $\mu$ on $(\Omega_{S},\mathcal{F}_{S})$,
and $f\in C(\Omega_{S})$.

Let $\theta$ be the measure on $(\Omega_{S},\mathcal{F}_{S})$ with
\[
\theta(E)=\int_{E}e^{P-\phi}\cdot\frac{h\circ T}{h}\: d\mathbb{P}
\]
for each $E\in\mathcal{F}_{S}$. Let $n\geq1$ and $a_{0}\cdot...\cdot a_{n-1}=w\in\mathcal{A}_{S}^{*}$,
then
\begin{multline}
\theta[w]=\int_{[w]}e^{P-\phi}\cdot h\circ T\: d\nu=\int_{[w]}e^{-\phi}\cdot h\circ T\: d(L_{\phi}^{*}\nu)=\\
=\int_{[w]}L_{\phi}(1_{[w]}\cdot e^{-\phi}\cdot h\circ T)\: d\nu=\int\sum_{y\in T^{-1}x}e^{\phi(y)}1_{[w]}(y)e^{-\phi(y)}h(Ty)\: d\nu(x)=\\
=\int1_{[w]}(a_{0}\cdot x)h(x)\: d\nu(x)=\int1_{[a_{1}\cdot...\cdot a_{n-1}]}(x)h(x)\: d\nu(x)=\mathbb{P}[a_{1}\cdot...\cdot a_{n-1}]=\mathbb{P}\circ T[w]\:.\label{E13}
\end{multline}

Set $\mathcal{L}=\{E\in\mathcal{F}_{S}\::\:\mathbb{P}\circ T(E)=\theta(E)\}$
and $\mathcal{P}=\{[w]\::\: w\in\mathcal{A}_{S}^{*}\}\cup\{\emptyset\}$.
From (\ref{E13}) it follows that $\Omega_{S}\in\mathcal{L}$, and
so it is easy to check that $\mathcal{L}$ is a $\lambda$-system.
Also, $\mathcal{P}$ is a $\pi$-system and from (\ref{E13}) we get
that $\mathcal{P}\subset\mathcal{L}$. From the $\pi-\lambda$ theorem
(see \cite{key-9}) it follows that $\mathcal{F}_{S}=\sigma(\mathcal{P})\subset\mathcal{L}$,
which shows that $\mathbb{P}\circ T=\theta$. This shows that $d\mathbb{P}=e^{\phi-P}\cdot\frac{h}{h\circ T}\: d\mathbb{P}\circ T$,
as required. $\square$

\emph{Proof of Corollary \ref{C2}:} From Theorem \ref{T4} and since
$T^{r}\omega=\omega$, it follows that
\[
\prod_{j=0}^{r-1}J(T^{j}\omega)=\prod_{j=0}^{r-1}e^{\phi(T^{j}\omega)-P}\frac{h(T^{j}\omega)}{h(T^{j+1}\omega)}=\exp(\sum_{j=0}^{r-1}(\phi(T^{j}\omega)-P))
\]
Now since $J=e^{\phi-P}\cdot\frac{h}{h\circ T}$ is a continuous function,
the claim follows from Corollary \ref{C1}. $\square$

We now turn to the proof of Theorem \ref{T5}. We shall need some
additional notations. For each $a\in\mathcal{A}$ set $K_{a}=f_{a}(I_{a})$
and let $g_{a}:K_{a}\rightarrow I_{a}$ be the inverse function of
$f_{a}$. From the assumptions on $f$ it follows that $g_{a}\in C^{1}(K_{a})$,
and $g_{a}$ is strictly monotonic. Let $a_{0}\cdot...\cdot a_{r-1}=w\in\mathcal{A}_{S}^{*}$,
then 
\[
g_{a_{j}}(K_{a_{j}})=I_{a_{j}}\subset f_{a_{j-1}}(I_{a_{j-1}})=K_{a_{j-1}}
\]
for each $1\leq j<r$, so the function $g_{w}:K_{a_{r-1}}\rightarrow I$
given by $g_{w}=g_{a_{0}}\circ...\circ g_{a_{r-1}}$ is well defined.
For such a $w$ we set $I_{w}=g_{w}(K_{a_{r-1}})$. 

The following lemmas will be needed in the proof of Theorem \ref{T5}.
\begin{lem*}
$\Theta(\omega)$ equals the unique element in $\bigcap_{n\geq1}I_{\omega_{0}\cdot...\cdot\omega_{n-1}}$
for each $\omega\in\Omega_{S}$, and
\begin{equation}
\Theta(a\cdot\omega)=g_{a}(\Theta(\omega))\mbox{ for each }a\in\mathcal{A}\mbox{ and }\omega\in T[a]\:.\label{E14}
\end{equation}

\end{lem*}

\emph{Proof:} First we shall prove by induction on $n$ that 
\begin{equation}
\bigcap_{j=0}^{n-1}f^{-j}(U_{a_{j}})=g_{a_{0}\cdot...\cdot a_{n-2}}(U_{a_{n-1}})\label{E32}
\end{equation}
for each $n\geq1$ and $a_{0}\cdot...\cdot a_{n-1}\in\mathcal{A}_{S}^{*}$.
This is clear for $n=1$ since both sides are equal to $U_{a_{0}}$.
Assume (\ref{E32}) holds for $n\geq1$ and let $a_{0}\cdot...\cdot a_{n}\in\mathcal{A}_{S}^{*}$,
then
\begin{multline*}
\bigcap_{j=0}^{n}f^{-j}(U_{a_{j}})=U_{a_{0}}\cap f^{-1}(\bigcap_{j=0}^{n-1}f^{-j}(U_{a_{j+1}}))=\\
=U_{a_{0}}\cap f^{-1}(g_{a_{1}\cdot...\cdot a_{n-1}}(U_{a_{n}}))=f_{a}^{-1}(g_{a_{1}\cdot...\cdot a_{n-1}}(U_{a_{n}}))=g_{a_{0}\cdot...\cdot a_{n-1}}(U_{a_{n}})
\end{multline*}
so (\ref{E32}) holds for all $n\geq1$.

From this it follows that for each $\omega\in\Omega_{S}$
\[
\bigcap_{n\geq1}\overline{\bigcap_{j=0}^{n-1}f^{-j}(U_{\omega_{j}})}=\bigcap_{n\geq1}\overline{g_{\omega_{0}\cdot...\cdot\omega_{n-2}}(U_{\omega_{n-1}})}=\bigcap_{n\geq1}I_{\omega_{0}\cdot...\cdot\omega_{n-1}}
\]
and so, from the definition of $\Theta$ found in Proposition 1.2
in \cite{key-10}, $\Theta(\omega)$ equals the unique element in
$\bigcap_{n\geq1}I_{\omega_{0}\cdot...\cdot\omega_{n-1}}$. Let $a\in\mathcal{A}$,
$\omega\in T[a]$ and $\gamma=a\cdot\omega$, then since 

\[
\bigcap_{n\geq1}I_{\gamma_{0}\cdot...\cdot\gamma_{n-1}}=g_{a}(\bigcap_{n\geq1}I_{\omega_{0}\cdot...\cdot\omega_{n-1}})
\]
it follows that $\Theta(\gamma)=g_{a}(\Theta(\omega))$, and the lemma
is proved. $\square$

\begin{lem*}
Let $a_{0}\cdot...\cdot a_{n-1}=w\in\mathcal{A}_{S}^{*}$, then 
\begin{equation}
\Phi^{-1}[w]=\tilde{I}\cap I_{w}\label{E34}
\end{equation}

\end{lem*}

\emph{Proof:} For each $0\leq m\leq n$ set
\[
F_{m}=\{x\in\tilde{I}\::\: f^{j}(x)\in I_{a_{j}}\mbox{ for each }m\leq j<n\}
\]
and
\[
E_{m}=\{x\in\tilde{I}\::\: f_{a_{j-1}}\circ...\circ f_{a_{0}}(x)\in I_{a_{j}}\mbox{ for each }0\le j<m\}.
\]
Clearly $F_{m}\cap E_{m}=F_{m+1}\cap E_{m+1}$ for each $0\leq m<n$.
Also $\Phi^{-1}[w]=F_{0}\cap E_{0}$, since
\begin{multline*}
\Phi^{-1}[w]=\{x\in\tilde{I}\::\:\xi_{j}(x)=a_{j}\mbox{ for all }0\leq j<n\}=\\
=\{x\in\tilde{I}\::\: f^{j}(x)\in U_{a_{j}}\mbox{ for all }0\leq j<n\}=\\
=\{x\in\tilde{I}\::\: f^{j}(x)\in I_{a_{j}}\mbox{ for all }0\leq j<n\}
\end{multline*}
which shows that $\Phi^{-1}[w]=F_{n}\cap E_{n}$, and so 
\begin{multline*}
\Phi^{-1}[w]=\{x\in\tilde{I}\::\: f_{a_{j-1}}\circ...\circ f_{a_{0}}(x)\in I_{a_{j}}\mbox{ for each }0\le j<n\}=\\
=\{x\in\tilde{I}\::\: x\in g_{a_{0}}\circ...\circ g_{a_{j}}(K_{a_{j}})\mbox{ for all }0\le j<n\}.
\end{multline*}

For each $1\leq j<n$ it holds that $g_{a_{j}}(K_{a_{j}})\subset K_{a_{j-1}}$,
so 
\[
g_{a_{0}}\circ...\circ g_{a_{j}}(K_{a_{j}})\subset g_{a_{0}}\circ...\circ g_{a_{j-1}}(K_{a_{j-1}})
\]
and so 
\[
\Phi^{-1}[w]=\{x\in\tilde{I}\::\: x\in g_{a_{0}}\circ...\circ g_{a_{n-1}}(K_{a_{n-1}})\}=\tilde{I}\cap I_{w}
\]
as desired. $\square$

\begin{lem*}
For every $\mathbb{P}\circ T$-integrable function $h$ it holds that
\begin{equation}
\int h\: d\mathbb{P}\circ T=\sum_{a\in\mathcal{A}}\int_{T[a]}h(a\cdot\omega)\: d\mathbb{P}(\omega).\label{E33}
\end{equation}

\end{lem*}

\emph{Proof:} Let $E\in\mathcal{F}_{S}$, then
\begin{multline*}
\int1_{E}\: d\mathbb{P}\circ T=\sum_{a\in\mathcal{A}}\mathbb{P}(T(E\cap[a]))=\\
=\sum_{a\in\mathcal{A}}\int1_{T(E\cap[a])}(\omega)\: d\mathbb{P}(\omega)=\sum_{a\in\mathcal{A}}\int_{T[a]}1_{E}(a\cdot\omega)\: d\mathbb{P}(\omega)
\end{multline*}
hence the lemma holds for indicator functions. From linearity the
lemma follows for simple functions, from the monotone convergence
theorem it follows for positive functions, and by writing $h=(h\vee0)-(-h\vee0)$
the lemma follows for every $\mathbb{P}\circ T$-integrable function
$h$. $\square$

\emph{Proof of Theorem \ref{T5}:} For each $\omega\in\Omega_{S}$
set $h(\omega)=\frac{p(\Theta(\omega))}{f_{\omega_{0}}'(\Theta(\omega))p(\Theta(T\omega))}$.
Let $X=\Phi(\tilde{I})$, then $\Theta:X\rightarrow\tilde{I}$ is
a bijection whose inverse is $\Phi$, and also 
\[
\mathbb{P}(X)=\mu(\Phi^{-1}(\Phi(\tilde{I})))=\mu(\tilde{I})=1\:.
\]
Let $a_{0}\cdot...\cdot a_{r-1}=w\in\mathcal{A}_{S}^{*}$, and set
$\alpha=\int_{[w]}h(\omega)\: d\mathbb{P}\circ T$. From (\ref{E33}),
\begin{multline}
\alpha=\sum_{a\in\mathcal{A}}\int_{T[a]}1_{[w]}(a\cdot\omega)h(a\cdot\omega)\: d\mathbb{P}(\omega)=\\
=\int_{T[a_{0}]}1_{[w]}(a_{0}\cdot\omega)h(a_{0}\cdot\omega)\: d\mathbb{P}(\omega)\:.\label{E35}
\end{multline}
It holds that
\[
\Theta(T[a_{0}]\cap X)=\Theta\{\omega\in X\::\: S_{a_{0},\omega_{0}}=1\}=\tilde{I}\cap(\bigcup\{U_{b}\::\: S_{a_{0},b}=1\})=\tilde{I}\cap K_{a_{0}}
\]
so $\omega\in T[a_{0}]\cap X$ if and only if $\Theta(\omega)\in\tilde{I}\cap K_{a_{0}}$.
From (\ref{E34}) it follows that
\[
\Theta([w]\cap X)=\Phi^{-1}[w]=I_{w}\cap\tilde{I},
\]
so $a_{0}\cdot\omega\in[w]\cap X$ if and only if $\Theta(a_{0}\cdot\omega)\in I_{w}\cap\tilde{I}$
for each $\omega\in T[a_{0}]$. Also, if $\omega\in T[a_{0}]\cap X$
then $a_{0}\cdot\omega\in X$. From this, from $\mathbb{P}(X)=1$,
and from (\ref{E35}),
\[
\alpha=\int_{X}1_{K_{a_{0}}}(\Theta(\omega))\cdot1_{I_{w}}(\Theta(a_{0}\cdot\omega))\cdot\frac{p(\Theta(a_{0}\cdot\omega))}{f_{a_{0}}'(\Theta(a_{0}\cdot\omega))p(\Theta(\omega))}\: d\mathbb{P}(\omega).
\]
This together with (\ref{E14}) gives
\begin{multline*}
\alpha=\int_{X}1_{K_{a_{0}}}(\Theta(\omega))\cdot1_{I_{w}}(g_{a_{0}}(\Theta(\omega)))\cdot\frac{p(g_{a_{0}}(\Theta(\omega)))}{f_{a_{0}}'(g_{a_{0}}(\Theta(\omega)))p(\Theta(\omega))}\: d\mathbb{P}(\omega)=\\
=\int_{\tilde{I}}1_{K_{a_{0}}}(x)\cdot1_{I_{w}}(g_{a_{0}}(x))\cdot\frac{p(g_{a_{0}}(x))}{f_{a_{0}}'(g_{a_{0}}(x))p(x)}\: d\mathbb{P}\circ\Theta^{-1}(x)
\end{multline*}
and since $g_{a_{0}}^{-1}(I_{w})\subset K_{a_{0}}$,
\[
\alpha=\int_{\tilde{I}}1_{g_{a_{0}}^{-1}(I_{w})}(x)\cdot\frac{p(g_{a_{0}}(x))}{f_{a_{0}}'(g_{a_{0}}(x))p(x)}\: d\mathbb{P}\circ\Theta^{-1}(x).
\]
Since $\mathbb{P}=\mu\circ\Phi^{-1}$ and $d\mu=p\: dm$, this shows
that
\begin{multline*}
\alpha=\int_{\tilde{I}}1_{g_{a_{0}}^{-1}(I_{w})}(x)\cdot\frac{p(g_{a_{0}}(x))}{f_{a_{0}}'(g_{a_{0}}(x))p(x)}\: d\mu(x)=\\
=\int1_{g_{a_{0}}^{-1}(I_{w})}(x)\cdot\frac{p(g_{a_{0}}(x))}{f_{a_{0}}'(g_{a_{0}}(x))}\: dm(x),
\end{multline*}
and since $g_{a_{0}}'=\frac{1}{f_{a_{0}}'\circ g_{a_{0}}}$,
\begin{multline*}
\alpha=\int1_{g_{a_{0}}^{-1}(I_{w})}(x)\cdot p(g_{a_{0}}(x))g_{a_{0}}'(x)\: dm(x)=\\
=\int_{I_{w}}p(y)\: dm(y)=\mu([I_{w}]\cap\tilde{I})=\mu(\Phi^{-1}[w])=\mathbb{P}[w].
\end{multline*}
This holds for any $w\in\mathcal{A}_{S}^{*}$, hence from the $\pi-\lambda$
theorem (see the end of the proof of Theorem \ref{T4}) it follows
that $\frac{d\mathbb{P}}{d\mathbb{P}\circ T}=h$, as required. $\square$\emph{$\newline$$\newline$Proof
of Corollary \ref{C3}:} From Theorem \ref{T5} and since $T^{r}\omega=\omega$,
it follows that
\[
\prod_{j=0}^{r-1}J(T^{j}\omega)=\prod_{j=0}^{r-1}\frac{p(\Theta(T^{j}\omega))}{f_{\omega_{j}}'(\Theta(T^{j}\omega))p(\Theta(T^{j+1}\omega))}=\left(\prod_{j=0}^{r-1}f_{\omega_{j}}'(\Theta(T^{j}\omega))\right)^{-1}.
\]
Now since the map $\omega\rightarrow\frac{p(\Theta(\omega))}{f_{\omega_{0}}'(\Theta(\omega))p(\Theta(T\omega))}$
is continuous, the claim follows from Corollary \ref{C1}. $\square$

\section{\label{S6}Proof of Theorem \ref{T6}}

\emph{Proof of Assertion (a):} Let $\epsilon>0$, let $\Gamma>0$
be as in (\ref{E2}), and let $M\geq1$ be such that $ne^{-\Gamma n}<1$
for all $n\geq M$. Let $n\geq M$, then from Lemma 3.2 in \cite{key-1}
it follows that for each $n\leq k\le N_{n}^{\omega}$,
\[
\mathbb{P}\{X_{k}^{A_{n}^{\omega}}=1\}\leq(1+\psi_{0})^{\ell}(\mathbb{P}(A_{n}^{\omega}))^{\ell}.
\]
 Also, for each $1\leq k<n$ it follows from (\ref{E2}) that 
\[
\mathbb{P}\{X_{k}^{A_{n}^{\omega}}=1\}\leq\mathbb{P}(A_{n}^{\omega})\leq e^{-\Gamma n}.
\]
Hence,
\begin{equation}
E[S_{N_{n}^{\omega}}^{A_{n}^{\omega}}]=\underset{k=1}{\overset{N_{n}^{\omega}}{\sum}}\mathbb{P}\{X_{k}^{A_{n}^{\omega}}=1\}\leq ne^{-\Gamma n}+N_{n}^{\omega}(1+\psi_{0})^{\ell}(\mathbb{P}(A_{n}^{\omega}))^{\ell}\leq1+(1+\psi_{0})^{\ell}t.\label{E15}
\end{equation}
Let $b\geq\frac{1+(1+\psi_{0})^{\ell}t}{\epsilon}$ be sufficiently
large such that $\mu_{n}[b,\infty)\leq\epsilon$ for each $1\leq n<M$.
For each $n\geq M$ it follows from (\ref{E15}) that 
\[
\mu_{n}[b,\infty)=\mathbb{P}\{S_{N_{n}^{\omega}}^{A_{n}^{\omega}}\geq b\}\leq\frac{E[S_{N_{n}^{\omega}}^{A_{n}^{\omega}}]}{b}\leq\epsilon
\]
which shows that $\{\mu_{n}\}_{n=1}^{\infty}$ is tight.

\emph{Proof of Assertion (b):} For each $n\geq1$ let $W_{n}\sim Pois(t(1-\rho_{A_{n}^{\omega}}))$
(where $\rho_{A_{n}^{\omega}}$ is defined in section \ref{SS5}).
From Theorem 2.3 in \cite{key-1}, it follows that for each $n\geq1$
there exist a sequence of i.i.d. random variables $\eta_{n,1},\eta_{n,2},...$,
independent of $W_{n}$, such that $\mathbb{P}\{\eta_{n,1}\in\{1,...,[\frac{n}{r}]\}\}=1$
and for $Z_{n}=\sum_{j=1}^{W_{n}}\eta_{n,j}$,
\begin{equation}
\underset{n\rightarrow\infty}{\lim}d_{TV}(\mu_{n},\mathcal{L}(Z_{n}))=0\:.\label{E16}
\end{equation}
For each $n\geq1$ let $\nu_{n}\in\mathcal{M}(\mathbb{N})$ be the
distribution of $Z_{n}$. From (\ref{E16}) and since $\{\mu_{n}\}_{n=1}^{\infty}$
is tight, it follows that $\{\nu_{n}\}_{n=1}^{\infty}$ is also tight.$\newline$For
each $n\geq1$ let $\theta_{n}$ be the distribution of $\eta_{n,1},\eta_{n,2},...$.
We shall now show that $\{\theta_{n}\}_{n=1}^{\infty}$ is tight.
Assume by contradiction that $\{\theta_{n}\}_{n=1}^{\infty}$ is not
tight, then there exist $\epsilon>0$ such that for each $k\geq1$
there exist $n_{k}\geq1$ with $\theta_{n_{k}}[k,\infty)>\epsilon$.
From Theorem 2.3 in \cite{key-1} it follows that $\underset{n\geq1}{\sup}\:\rho_{A_{n}^{\omega}}<1$,
so there exists $\delta>0$ such that
\[
\underset{n\geq1}{\inf}\:\mathbb{P}\{W_{n}>0\}=\underset{n\geq1}{\inf}\:(1-\exp(-t(1-\rho_{A_{n}^{\omega}})))>\delta,
\]
and so for each $k\geq1$,
\begin{multline*}
\nu_{n_{k}}[k,\infty)=\mathbb{P}\{Z_{n_{k}}\geq k\}\geq\mathbb{P}(\{W_{n_{k}}>0\}\cap\{\eta_{n_{k},1}\geq k\})=\\
=\mathbb{P}\{W_{n_{k}}>0\}\cdot\theta_{n_{k}}[k,\infty)>\delta\cdot\epsilon
\end{multline*}
which is a contradiction to the tightness of $\{\nu_{n}\}_{n=1}^{\infty}$,
and so $\{\theta_{n}\}_{n=1}^{\infty}$ must be tight.$\newline$Let
$\mu$ be a probability distribution on $\mathbb{R}$ such that $\mu_{n_{k}}\overset{d}{\Longrightarrow}\mu$
as $k\rightarrow\infty$, for some increasing sequence $\{n_{k}\}_{k=1}^{\infty}\subset\mathbb{N}_{+}$.
From (\ref{E16}) it follows that also $\nu_{n_{k}}\overset{d}{\Longrightarrow}\mu$
as $k\rightarrow\infty$. For each $k\geq1$ set $\tau_{k}=t(1-\rho_{A_{n_{k}}^{\omega}})$,
then $W_{n_{k}}\sim Pois(\tau_{k})$ and
\[
0<t(1-\underset{n\geq1}{\sup}\:\rho_{A_{n}^{\omega}})\leq\tau_{k}\leq t
\]
for each $k\geq1$. From this it follows, by moving to a sub-sequence
without changing notation, that we can assume $\tau_{k}\overset{k}{\rightarrow}\tau$
for some $0<\tau\le t$. Also, since $\{\theta_{n}\}_{n=1}^{\infty}$
is tight and from Theorem 25.10 in \cite{key-9}, by moving to a sub-sequence
without changing notation, we can assume $\theta_{n_{k}}\overset{d}{\Longrightarrow}\theta$
for some probability distribution $\theta$ on $\mathbb{R}$. Since
$\theta_{n_{k}}\in\mathcal{M}(\mathbb{N})$ for each $k\geq1$, it
follows that $\theta\in\mathcal{M}(\mathbb{N})$.$\newline$It follows
from Theorem 26.3 in \cite{key-9} and from (\ref{E8}), that for
each $x\in\mathbb{R}$
\begin{multline*}
\varphi_{\mu}(x)=\underset{k\rightarrow\infty}{\lim}\varphi_{\nu_{n_{k}}}(x)=\underset{k\rightarrow\infty}{\lim}\exp(\tau_{k}\cdot(\varphi_{\theta_{n_{k}}}(x)-1))=\\
=\exp(\tau\cdot(\varphi_{\theta}(x)-1))=\varphi_{CP(\tau,\theta)}(x)
\end{multline*}
which shows that $\mu=CP(\tau,\theta)$, and completes the proof.
$\square$

\section{\label{S7}Proof of Theorem \ref{T7}}

\emph{Proof of Assertion (a):} Let $g_{0},...,g_{n-1}\in G$, then
since $\mathbb{P}$ is $T$-invariant,
\begin{multline*}
\mathbb{P}_{0}(T^{-1}[g_{0},...,g_{n-1}])=\mathbb{P}(\Phi^{-1}\{\omega\::\:\omega_{j+1}=g_{j}\mbox{ for each }0\le j<n\})=\\
=\mathbb{P}(\bigcap_{j=0}^{n-1}\{\omega\::\:\omega_{j+1}+...+\omega_{j+N}=g_{j}\})=\mathbb{P}(T^{-1}(\bigcap_{j=0}^{n-1}\{\omega\::\:\omega_{j}+...+\omega_{j+N-1}=g_{j}\}))=\\
=\mathbb{P}(\bigcap_{j=0}^{n-1}\{\omega\::\:\omega_{j}+...+\omega_{j+N-1}=g_{j}\})=\mathbb{P}(\Phi^{-1}[g_{0},...,g_{n-1}])=\mathbb{P}_{0}[g_{0},...,g_{n-1}].
\end{multline*}
Since $\mathcal{F}$ is generated by the $\pi$-system of all cylinders,
it follows from the $\pi-\lambda$ theorem, that $\mathbb{P}_{0}$
is $T$-invariant. $\square$

For the proof of Assertion (b) the following lemma will be needed.
\begin{lem}
\label{L4}Let $l\in\mathbb{N}$ and $C\geq0$ be such that
\begin{equation}
|\mathbb{P}_{0}(A\cap T^{-n-l}B)-\mathbb{P}_{0}(A)\mathbb{P}_{0}(B)|\leq C\mathbb{P}_{0}(A)\mathbb{P}_{0}(B)\label{E19}
\end{equation}
for each $n$-cylinder $A\in\mathcal{F}_{\{0,...,n-1\}}$ and cylinder
$B\in\mathcal{F}$. Then
\[
|\mathbb{P}_{0}(E\cap T^{-n-l}F)-\mathbb{P}_{0}(E)\mathbb{P}_{0}(F)|\leq C\mathbb{P}_{0}(E)\mathbb{P}_{0}(F)
\]
for each $E\in\mathcal{F}_{\{0,...,n-1\}}$ and $F\in\mathcal{F}$.
\end{lem}

\emph{Proof: }Since $\mathcal{F}$ is generated by cylinders the result
follows. $\square$

\emph{Proof of Assertion (b): }Let $n,k\geq1$, $g_{0},...,g_{n+k-1}\in G$,
$A=[g_{0},...,g_{n-1}]$, $B=[g_{n},...,g_{n+k-1}]$ and $\lambda=\min\{p_{g}\::\: g\in G\}$.
We shall first show by induction on $m$ that for each $0\leq m<N$
\begin{equation}
\mathbb{P}_{0}[g_{0},...,g_{n-N}]\leq\lambda^{-m}\mathbb{P}_{0}[g_{0},...,g_{n+m-N}]\:.\label{E17}
\end{equation}
For $m=0$, (\ref{E17}) is obvious. Assume (\ref{E17}) is true for
some $0\leq m<N-1$. If $n+m<N-1$ then $[g_{0},...,g_{n+m+1-N}]=\Omega$
and (\ref{E17}) is obvious, hence we can assume that $n+m\geq N-1$.
Set
\[
\mathcal{E}=\{R\in\mathcal{F}\::\: R=[x_{0},...,x_{n+m-1}]\mbox{ and }x_{j}+...+x_{j+N-1}=g_{j}\mbox{ for }0\leq j\leq n+m-N\}
\]
then
\[
\Phi^{-1}[g_{0},...,g_{n+m-N}]=\bigcup_{R\in\mathcal{E}}R\:.
\]
For each $[x_{0},...,x_{n+m-1}]=R\in\mathcal{E}$ set
\[
y_{R}=g_{n+m+1-N}-x_{n+m-1}-...-x_{n+m+1-N}
\]
and $Q_{R}=[x_{0},...,x_{n+m-1},y_{R}]$. Then $Q_{R}\subset\Phi^{-1}[g_{0},...,g_{n+m+1-N}]$,
and so
\begin{multline*}
\mathbb{P}_{0}[g_{0},...,g_{n-N}]\overset{\underbrace{\mbox{i.h.}}}{\leq}\lambda^{-m}\mathbb{P}_{0}[g_{0},...,g_{n+m-N}]=\lambda^{-m}\mathbb{P}(\Phi^{-1}[g_{0},...,g_{n+m-N}])=\\
=\lambda^{-m}\sum_{R\in\mathcal{E}}\mathbb{P}(R)\leq\lambda^{-m-1}\sum_{R\in\mathcal{E}}\mathbb{P}(Q_{R})=\lambda^{-m-1}\mathbb{P}(\bigcup_{R\in\mathcal{E}}Q_{R})\leq\\
\leq\lambda^{-m-1}\mathbb{P}(\Phi^{-1}[g_{0},...,g_{n+m+1-N}])=\lambda^{-m-1}\mathbb{P}_{0}[g_{0},...,g_{n+m+1-N}]
\end{multline*}
and the induction is complete.

Let $0\leq l<N-1$, then from (\ref{E17}) with $m=N-1$, 
\begin{multline*}
|\mathbb{P}_{0}(A\cap T^{-n-l}B)-\mathbb{P}_{0}(A)\mathbb{P}_{0}(B)|\leq\mathbb{P}_{0}([g_{0},...,g_{n-N}]\cap T^{-n-l}B)+\mathbb{P}_{0}(A)\mathbb{P}_{0}(B)=\\
=\mathbb{P}_{0}[g_{0},...,g_{n-N}]\mathbb{P}_{0}(B)+\mathbb{P}_{0}(A)\mathbb{P}_{0}(B)\leq(1+\lambda^{-N+1})\mathbb{P}_{0}(A)\mathbb{P}_{0}(B)
\end{multline*}
and for $l\geq N-1$,
\[
\mathbb{P}_{0}(A\cap T^{-n-l}B)=\mathbb{P}_{0}(A)\mathbb{P}_{0}(B)\:.
\]
This together with Lemma \ref{L4} shows that \noun{$\mathbb{P}_{0}$
}is $\psi$-mixing. $\square$

\emph{Proof of Assertion (c):} Let $h\in G$ be such that $p_{h}>p_{g}$
for all $g\in G\setminus\{h\}$. The following notation will be needed.
For $m\in\mathbb{N}$ and $g\in G$ set 
\[
m\cdot g:=\underset{m\mbox{ times}}{\underbrace{g+...+g}}
\]
Let $G^{*}$ be the set of finite words over $G$ (when $G$ is thought
of as an alphabet). As before, for each $u,w\in G^{*}$ and $k\geq0$
let $u\cdot w\in G^{*}$ be the concatenation of $u$ and $w$, and
let $w^{k}\in G^{*}$ be the concatenation of $w$ with itself $k$
times.

Let $r\in G\setminus\{h\}$ be such that $p_{r}\geq p_{g}$ for all
$g\in G\setminus\{h\}$, and set $s=(N-1)\cdot h+r$. We will prove
(c) by showing that the limit $\underset{n}{\lim}\:\mathbb{P}_{0}\{[s^{n+1}]\mid[s^{n}]\}$
does not exist.$\newline$Set
\[
\mathcal{E}=\{R\in\mathcal{F}\::\: R=[g_{0},...,g_{N-1}]\mbox{ and }g_{0}+...+g_{N-1}=s\}
\]
and for each $0\leq j<N$ set $H_{j}=[h^{j}\cdot r\cdot h^{N-1-j}]$,
then since $G$ is abelian $H_{0},...,H_{N-1}\in\mathcal{E}$. For
each $n\geq1$, let $a_{n},b_{n}\in\mathbb{N}$ be such that $n+N-1=a_{n}\cdot N+b_{n}$
and $0\leq b_{n}<N$.$\newline$We shall now show by induction on
$n\geq1$ that $\Phi^{-1}[s^{n}]=\bigcup_{R\in\mathcal{E}}R^{(n+N-1)/N}$.
For $n=1$ this follows directly from the definition of $\mathcal{E}$.
Let $n\geq1$ and assume we know that $\Phi^{-1}[s^{n}]=\bigcup_{R\in\mathcal{E}}R^{(n+N-1)/N}$,
then
\begin{multline}
\Phi^{-1}[s^{n+1}]=\Phi^{-1}[s^{n}]\cap\Phi^{-1}\{\omega_{n}=s\}=\\
=(\bigcup_{R\in\mathcal{E}}R^{(n+N-1)/N})\cap\Phi^{-1}\{\omega_{n}=s\}=\bigcup_{R\in\mathcal{E}}(R^{(n+N-1)/N}\cap\Phi^{-1}\{\omega_{n}=s\})\:.\label{E21}
\end{multline}
Let $[g_{0},...,g_{N-1}]=R\in\mathcal{E}$, then since $G$ is abelian
it follows for $\omega\in R^{(n+N-1)/N}$ that
\[
\omega_{n-1}+...+\omega_{n+N-2}=g_{b_{n}}+...+g_{N-1}+g_{0}+...+g_{b_{n}-1}=s
\]
and for $\omega\in\Phi^{-1}\{\omega_{n}=s\}$,
\[
\omega_{n}+...+\omega_{n+N-1}=s.
\]
Hence, for $\omega\in R^{(n+N-1)/N}\cap\Phi^{-1}\{\omega_{n}=s\}$,
\[
\omega_{n-1}+...+\omega_{n+N-2}=s=\omega_{n}+...+\omega_{n+N-1}
\]
which shows that $\omega_{n-1}=\omega_{n+N-1}$, and so $\omega\in R^{(n+N)/N}$.
This shows that $R^{(n+N-1)/N}\cap\Phi^{-1}\{\omega_{n}=s\}\subset R^{(n+N)/N}$.
On the other hand, if $\omega\in R^{(n+N)/N}$ then
\[
\omega_{n}+...+\omega_{n+N-1}=g_{b_{n+1}}+...+g_{N-1}+g_{0}+...+g_{b_{n+1}-1}=s
\]
so $\omega\in R^{(n+N-1)/N}\cap\Phi^{-1}\{\omega_{n}=s\}$, which
shows that 
\[
R^{(n+N-1)/N}\cap\Phi^{-1}\{\omega_{n}=s\}=R^{(n+N)/N}.
\]
Now from (\ref{E21}) it follows that $\Phi^{-1}[s^{n+1}]=\bigcup_{R\in\mathcal{E}}R^{(n+N)/N}$
and the induction is complete.$\newline$Let $n\geq1$, then
\begin{multline}
\mathbb{P}_{0}[s^{n}]=\mathbb{P}(\Phi^{-1}[s^{n}])=\mathbb{P}(\bigcup_{R\in\mathcal{E}}R^{(n+N-1)/N})=\sum_{R\in\mathcal{E}}\mathbb{P}(R^{(n+N-1)/N})=\\
=\sum_{[g_{0},...,g_{N-1}]\in\mathcal{E}}\:(p_{g_{0}}\cdot...\cdot p_{g_{N-1}})^{a_{n}}\cdot p_{g_{0}}\cdot...\cdot p_{g_{b_{n}-1}}\label{E22}
\end{multline}
Set $\mathcal{Q}=\{H_{0},...,H_{N-1}\}$ and let $[g_{0},...,g_{N-1}]\in\mathcal{E}\setminus\mathcal{Q}$,
we shall now show that $\frac{p_{g_{0}}\cdot...\cdot p_{g_{N-1}}}{p_{h}^{N-1}\cdot p_{r}}<1$.
Let $t\in G\setminus\{r\}$, then
\[
(N-1)h+t\ne(N-1)h+r=s
\]
and so $[g_{0},...,g_{N-1}]\ne[h^{j}\cdot t\cdot h^{N-1-j}]$ for
each $0\leq j<N$. Since $[g_{0},...,g_{N-1}]\notin\mathcal{Q}$,
it follows that there exist $0\leq i<j<N$ with $g_{i},g_{j}\ne h$.
Since $p_{h}>p_{g}$ and $p_{r}\geq p_{g}$ for each $g\in G\setminus\{h\}$,
it follows that 
\begin{equation}
\frac{p_{g_{0}}\cdot...\cdot p_{g_{N-1}}}{p_{h}^{N-1}\cdot p_{r}}\leq\frac{p_{g_{i}}\cdot p_{g_{j}}}{p_{h}\cdot p_{r}}<1.\label{E23}
\end{equation}
From $\mathcal{Q}\subset\mathcal{E}$, (\ref{E22}) and (\ref{E23})
we get that
\begin{multline*}
1\leq\underset{n}{\limsup}\:\frac{\mathbb{P}_{0}[s^{n}]}{\sum_{R\in\mathcal{Q}}\mathbb{P}(R^{(n+N-1)/N})}\leq\\
\leq1+\underset{n}{\limsup}\:\sum_{[g_{0},...,g_{N-1}]\in\mathcal{E}\setminus\mathcal{Q}}\frac{(p_{g_{0}}\cdot...\cdot p_{g_{N-1}})^{a_{n}}}{(p_{h}^{N-1}\cdot p_{r})^{a_{n}+1}}=1,
\end{multline*}
and so
\begin{equation}
\underset{n}{\lim}\:\mathbb{P}_{0}\{[s^{n+1}]\mid[s^{n}]\}\cdot(\frac{\sum_{R\in\mathcal{Q}}\mathbb{P}(R^{(n+N)/N})}{\sum_{R\in\mathcal{Q}}\mathbb{P}(R^{(n+N-1)/N})})^{-1}=1.\label{E24}
\end{equation}
Let $n\geq1$, then if $b_{n}=0$,
\begin{multline}
\frac{\sum_{R\in\mathcal{Q}}\mathbb{P}(R^{(n+N)/N})}{\sum_{R\in\mathcal{Q}}\mathbb{P}(R^{(n+N-1)/N})}=\\
=\frac{(p_{h}^{N-1}\cdot p_{r})^{a_{n}}\cdot p_{r}+(N-1)(p_{h}^{N-1}\cdot p_{r})^{a_{n}}\cdot p_{h}}{N\cdot(p_{h}^{N-1}\cdot p_{r})^{a_{n}}}=\frac{p_{r}+(N-1)p_{h}}{N}\label{E25}
\end{multline}
and if $b_{n}=N-1$ then
\begin{multline}
\frac{\sum_{R\in\mathcal{Q}}\mathbb{P}(R^{(n+N)/N})}{\sum_{R\in\mathcal{Q}}\mathbb{P}(R^{(n+N-1)/N})}=\\
=\frac{N\cdot(p_{h}^{N-1}\cdot p_{r})^{a_{n}+1}}{(N-1)(p_{h}^{N-1}\cdot p_{r})^{a_{n}}\cdot p_{h}^{N-2}\cdot p_{r}+(p_{h}^{N-1}\cdot p_{r})^{a_{n}}\cdot p_{h}^{N-1}}=\frac{N\cdot p_{h}\cdot p_{r}}{(N-1)\cdot p_{r}+p_{h}}.\label{E26}
\end{multline}
Now if
\[
\frac{p_{r}+(N-1)p_{h}}{N}=\frac{N\cdot p_{h}\cdot p_{r}}{(N-1)\cdot p_{r}+p_{h}},
\]
a direct computation shows that $p_{h}=p_{r}$, so according to our
assumptions it must hold that
\begin{equation}
\frac{p_{r}+(N-1)p_{h}}{N}\ne\frac{N\cdot p_{h}\cdot p_{r}}{(N-1)\cdot p_{r}+p_{h}}.\label{E27}
\end{equation}
From (\ref{E24}), (\ref{E25}), (\ref{E26}) and (\ref{E27}) it
follows that $\underset{n}{\lim}\:\mathbb{P}_{0}\{[s^{n+1}]\mid[s^{n}]\}$
does not exist, and the theorem is proved. $\square$

\section{\label{S8}Proof of Theorem \ref{T8}}

For $b_{0},...,b_{n-1}\in\{0,1\}$ we write
\[
[b_{0},...,b_{n-1}]_{0}=\{\omega\in\Omega_{0}\::\:\omega_{j}=b_{j}\mbox{ for each }0\leq j<n\}.
\]

\emph{Proof of Assertion (a):} Let $b_{0},...,b_{n-1}\in\{0,1\}$,
then since $\mathbb{P}$ is $T$-invariant,
\begin{multline*}
\mathbb{P}_{0}(T_{0}^{-1}[b_{0},...,b_{n-1}]_{0})=\mathbb{P}(\Theta^{-1}\{\omega\in\Omega_{0}\::\:\omega_{j+1}=b_{j}\mbox{ for each }0\leq j<n\})=\\
=\mathbb{P}\{\omega\in\Omega\::\:\theta(\omega_{j+1},\omega_{j+2})=b_{j}\mbox{ for each }0\leq j<n\}=\\
=\mathbb{P}(T^{-1}\{\omega\in\Omega\::\:\theta(\omega_{j},\omega_{j+1})=b_{j}\mbox{ for each }0\leq j<n\})=\\
=\mathbb{P}\{\omega\in\Omega\::\:\theta(\omega_{j},\omega_{j+1})=b_{j}\mbox{ for each }0\leq j<n\}=\\
=\mathbb{P}(\Theta^{-1}\{\omega\in\Omega_{0}\::\:\omega_{j}=b_{j}\mbox{ for each }0\leq j<n\})=\mathbb{P}_{0}[b_{0},...,b_{n-1}]_{0}.
\end{multline*}
Since $\mathcal{F}_{0}$ is generated by the $\pi$-system of all
cylinders, it follows that $\mathbb{P}_{0}$ is $T_{0}$-invariant.
$\square$

For the proof of Assertion (b) we shall need the following lemma.
\begin{lem}
\label{L5}Let $b_{0},...,b_{n-1}\in\{0,1\}$, $D=\Theta^{-1}[b_{0},...,b_{n-1}]_{0}$
and integers $1\leq r<k\leq s$ be given, then 
\[
\mathbb{P}(\{\omega_{0}=s\}\cap D)\leq\mathbb{P}(\{\omega_{0}=r\}\cap D)+\mathbb{P}(\{\omega_{0}=k\}\cap D).
\]

\end{lem}

\emph{Proof:} If $b_{j}=1$ for each $0\leq j<n$ then
\begin{multline*}
\mathbb{P}(\{\omega_{0}=s\}\cap D)=\mathbb{P}\{\omega_{j}=s+j\::\:\mbox{for each }0\leq j\leq n\}=\prod_{j=0}^{n}2^{-s-j}\leq\\
\leq\prod_{j=0}^{n}2^{-k-j}=\mathbb{P}\{\omega_{j}=k+j\::\:\mbox{for each }0\leq j\leq n\}=\mathbb{P}(\{\omega_{0}=k\}\cap D).
\end{multline*}
Hence, we can assume that there exist $0\leq j_{0}<n$ such that $b_{j_{0}}=0$
and $b_{j}=1$ for each $0\leq j<j_{0}$. Set
\[
D'=\{\omega\in\Omega\::\:\theta(\omega_{j},\omega_{j+1})=b_{j}\mbox{ for each }j_{0}+1\le j<n\}.
\]
Then it follows that
\begin{multline*}
\mathbb{P}(\{\omega_{0}=r\}\cap D)+\mathbb{P}(\{\omega_{0}=k\}\cap D)=\\
=\mathbb{P}(\{\omega_{j}=r+j\::\:\mbox{for each }0\leq j\leq j_{0}\}\cap\{\omega_{j_{0}+1}\ne r+j_{0}+1\}\cap D')+\\
+\mathbb{P}(\{\omega_{j}=k+j\::\:\mbox{for each }0\leq j\leq j_{0}\}\cap\{\omega_{j_{0}+1}\ne k+j_{0}+1\}\cap D')\geq\\
\geq\mathbb{P}\{\omega_{j}=s+j\::\:\mbox{for each }0\leq j\leq j_{0}\}(\mathbb{P}(\{\omega_{j_{0}+1}\ne r+j_{0}+1\}\cap D')+\mathbb{P}(\{\omega_{j_{0}+1}\ne k+j_{0}+1\}\cap D'))\geq\\
\geq\mathbb{P}\{\omega_{j}=s+j\::\:\mbox{for each }0\leq j\leq j_{0}\}\mathbb{P}(D')\geq\mathbb{P}(\{\omega_{0}=s\}\cap D)
\end{multline*}
and the lemma is proved. $\square$\emph{$\newline$$\newline$Proof
of Assertion (b):} Let $n\geq0$, $l\geq1$, and $b_{0},...,b_{n+l}\in\{0,1\}$.
Set 
\[
B_{1}=\{\omega\in\Omega\::\:\theta(\omega_{j},\omega_{j+1})=b_{j}\mbox{ for each }j\in\{0,...,n-1\}\}
\]
and
\[
B_{2}=\{\omega\in\Omega\::\:\theta(\omega_{j},\omega_{j+1})=b_{j}\mbox{ for each }j\in\{n+1,...,n+l\}\}.
\]
Then $B_{1}$ and $B_{2}$ are independent events, and so
\begin{multline}
\mathbb{P}_{0}([b_{0},...,b_{n-1},0]_{0}\cap T_{0}^{-(n+1)}[b_{n+1},...,b_{n+l}]_{0})\leq\mathbb{P}(B_{1}\cap B_{2})=\\
=\mathbb{P}(B_{1})\mathbb{P}(B_{2})=\mathbb{P}(B_{2})\sum_{k=1}^{\infty}\mathbb{P}(B_{1}\cap\{\omega_{n}=k\})\mathbb{P}\{\omega_{n+1}\ne k+1\}(\mathbb{P}\{\omega_{n+1}\ne k+1\})^{-1}\leq\\
\leq\mathbb{P}(B_{2})\sum_{k=1}^{\infty}\mathbb{P}(B_{1}\cap\{\omega_{n}=k\}\cap\{\omega_{n+1}\ne k+1\})\cdot2=2\mathbb{P}(B_{2})\mathbb{P}(B_{1}\cap\{\theta(\omega_{n},\omega_{n+1})=0\})=\\
=2\mathbb{P}_{0}[b_{0},...,b_{n-1},0]_{0}\mathbb{P}_{0}[b_{n+1},...,b_{n+l}]_{0}\:.\label{E28}
\end{multline}
In a similar manner it can be shown that
\begin{equation}
\mathbb{P}_{0}([b_{0},...,b_{n}]_{0}\cap T_{0}^{-(n+1)}[0,b_{n+2},...,b_{n+l}]_{0})\leq2\mathbb{P}_{0}[b_{0},...,b_{n}]_{0}\mathbb{P}_{0}[0,b_{n+2},...,b_{n+l}]_{0}\:.\label{E29}
\end{equation}
Now set
\[
B_{3}=\{\omega\in\Omega\::\:\theta(\omega_{j},\omega_{j+1})=b_{j}\mbox{ for each }j\in\{n+2,...,n+l\}\}.
\]
Then
\begin{multline}
\mathbb{P}_{0}([b_{0},...,b_{n-1},1]_{0}\cap T_{0}^{-(n+1)}[1,b_{n+2},...,b_{n+l}]_{0})=\\
=\sum_{k=1}^{\infty}\mathbb{P}(B_{1}\cap\{\omega_{n}=k\}\cap\{\omega_{n+1}=k+1\}\cap\{\omega_{n+2}=k+2\}\cap B_{3})=\\
=\sum_{k=1}^{\infty}\mathbb{P}(B_{1}\cap\{\omega_{n}=k\}\cap\{\omega_{n+1}=k+1\})\mathbb{P}(\{\omega_{n+2}=k+2\}\cap B_{3})\overset{\underbrace{\mbox{lemma }\ref{L5}}}{\leq}\\
\leq(\mathbb{P}(\{\omega_{n+2}=3\}\cap B_{3})+\mathbb{P}(\{\omega_{n+2}=2\}\cap B_{3}))\sum_{k=1}^{\infty}\mathbb{P}(B_{1}\cap\{\omega_{n}=k\}\cap\{\omega_{n+1}=k+1\})=\\
=(4\mathbb{P}(\{\omega_{n+1}=2\}\cap\{\omega_{n+2}=3\}\cap B_{3})+2\mathbb{P}(\{\omega_{n+1}=1\}\cap\{\omega_{n+2}=2\}\cap B_{3}))\mathbb{P}_{0}[b_{0},...,b_{n-1},1]_{0}\leq\\
\leq6\mathbb{P}_{0}[b_{0},...,b_{n-1},1]_{0}\mathbb{P}_{0}[1,b_{n+2}...,b_{n+l}]_{0}\quad.\label{E30}
\end{multline}
From (\ref{E28}), (\ref{E29}), (\ref{E30}) and Lemma \ref{L4}
it follows that
\[
|\mathbb{P}_{0}(E\cap T^{-n-1}F)-\mathbb{P}_{0}(E)\mathbb{P}_{0}(F)|\leq7\mathbb{P}_{0}(E)\mathbb{P}_{0}(F)
\]
for each $E\in\mathcal{F}_{0,\{0,...,n\}}$ and $F\in\mathcal{F}_{0}$.$\newline$Let
$s>n+1$, then
\[
\mathbb{P}_{0}([b_{0},...,b_{n}]_{0}\cap T_{0}^{-s}[b_{n+1},...,b_{n+l}]_{0})=\mathbb{P}_{0}[b_{0},...,b_{n}]_{0}\mathbb{P}_{0}[b_{n+1},...,b_{n+l}]_{0}
\]
so by Lemma \ref{L4},
\[
\mathbb{P}_{0}(E\cap T^{-s}F)=\mathbb{P}_{0}(E)\mathbb{P}_{0}(F)
\]
for each $E\in\mathcal{F}_{0,\{0,...,n\}}$ and $F\in\mathcal{F}_{0}$.
This shows that $\mathbb{P}_{0}$ is $\psi$-mixing. $\square$\emph{$\newline$$\newline$Proof
of Assertion (c):} For each integer $l\in\mathbb{N}$ set $f(l)=2^{-l}$.
Let $n\geq1$, then
\begin{multline*}
\mathbb{P}_{0}[1^{n}]_{0}=\sum_{k=1}^{\infty}\mathbb{P}[k,...,k+n]=\sum_{k=1}^{\infty}\prod_{j=k}^{k+n}2^{-j}=\\
=\sum_{k=1}^{\infty}f(\sum_{j=k}^{k+n}j)=\sum_{k=1}^{\infty}f(\frac{(n+1)(2k+n)}{2})=\\
=f(\frac{(n+1)n}{2})\sum_{k=1}^{\infty}(2^{-(n+1)})^{k}=f(\frac{(n+1)n}{2})\cdot2^{-(n+1)}\cdot\frac{1}{1-2^{-(n+1)}}\;.
\end{multline*}
Hence,
\begin{multline*}
\mathbb{P}_{0}\{[1^{n+1}]\mid[1^{n}]\}=\frac{f(\frac{(n+2)(n+1)}{2})\cdot2^{-(n+2)}\cdot\frac{1}{1-2^{-(n+2)}}}{f(\frac{(n+1)n}{2})\cdot2^{-(n+1)}\cdot\frac{1}{1-2^{-(n+1)}}}\leq\\
\leq f(\frac{(n+2)(n+1)}{2}-\frac{(n+1)n}{2})=2^{-(n+1)}\overset{n\rightarrow\infty}{\longrightarrow}0
\end{multline*}
and the theorem is proved. $\square$

\end{document}